\documentclass{aims_hyp}
\usepackage{amsmath,cite}
  \usepackage{paralist}
\usepackage{graphicx}  \usepackage{epstopdf}
 \usepackage[colorlinks=true]{hyperref}
\hypersetup{urlcolor=blue, citecolor=red}

\def\currentvolume{X}
 \def\currentissue{X}
  \def\currentyear{200X}
   \def\currentmonth{XX}
    \def\ppages{X--XX}

\newtheorem{theorem}{Theorem}[section]

\newtheorem{lemma}{Lemma}[section]
\newtheorem{proposition}{Proposition}[section]

\newtheorem{problem}{Problem}[section]
\theoremstyle{definition}
\newtheorem{definition}{Definition}[section]
\newtheorem{remark}{Remark}[section]

\numberwithin{figure}{section}

\allowdisplaybreaks

\usepackage{enumerate} \usepackage{array}

\def\bbm{\begin{bmatrix}}
	\def\ebm{\end{bmatrix}}
\def\bpm{\begin{pmatrix}}
	\def\epm{\end{pmatrix}}
\def\bvm{\begin{vmatrix}}
	\def\evm{\end{vmatrix}}
\def\bin{\begin{enumerate}}
	\def\ein{\end{enumerate}}
\def\bit{\begin{itemize}}
	\def\eit{\end{itemize}}
\def\bid{\begin{description}}
	\def\eid{\end{description}}

\renewcommand{\div}{\text{\sl div}}

\newcommand{\mr}{\mathbb{R}}
\newcommand{\mR}{\mathbb{R}}

\newcommand{\deq}{:=}
\newcommand{\bt}{\boldsymbol{\tau}}
\newcommand{\bn}{\boldsymbol{\nu}}

\newcommand{\bx}{\boldsymbol{\xi}}

\newcommand{\p}{\partial}
\newcommand{\gr}{\nabla}

\renewcommand{\r}{\rho}
\newcommand{\g}{\gamma}

\newcommand{\z}{\zeta}
\renewcommand{\O}{\Omega}

\renewcommand{\L}{\Lambda}

\newcommand{\Gw}{\Gamma_{\text{\rm wedge}}}

\newcommand{\Gso}{\Gamma_{\text{\rm sonic}}}

\newcommand{\Gsh}{\Gamma_{\text{\rm shock}}}
\newcommand{\hGsh}{\hat\Gamma_{\text{\rm shock}}}
\newcommand{\Gsy}{\Gamma_{\text{\rm sym}}}

\title[THE SHOCK REFLECTION-DIFFRACTION PROBLEM $\,\,$] 
{Uniqueness and Stability \\ for the Shock Reflection-Diffraction Problem\\ for Potential Flow}

\author[Gui-Qiang G. Chen, Mikhail Feldman and Wei Xiang]{}


 \email{chengq@maths.ox.ac.uk}
 \email{feldman@math.wisc.edu}
 \email{weixiang@cityu.edu.hk}

\thanks{The research of Gui-Qiang G. Chen was supported in part by the UK	Engineering and Physical Sciences Research Council Award
	EP/E035027/1 and
	EP/L015811/1, and the Royal Society--Wolfson Research Merit Award (UK).
	The research
	of Mikhail Feldman was supported in part by the National Science Foundation under
	Grants DMS-1764278 and DMS-1401490, and the Van Vleck Professorship Research Award
	by the University of Wisconsin-Madison.
	The research of Wei Xiang was supported in part by the
	Research Grants Council of the HKSAR, China (Project No. CityU 21305215,
	Project No. CityU 11332916, Project No. CityU 11304817, and Project No. CityU 11303518).}

\begin{document}
\maketitle

\centerline{\scshape Gui-Qiang G. Chen}
\medskip
{\footnotesize
 \centerline{Mathematical Institute, University of Oxford,
 	Oxford, OX2 6GG, UK}
}

\medskip

\centerline{\scshape Mikhail Feldman}
\medskip
{\footnotesize
 \centerline{Department of Mathematics, University of Wisconsin, Madison, WI 53706-1388, USA}
}

\medskip

\centerline{\scshape Wei Xiang}
\medskip
{\footnotesize
	\centerline{City University of Hong Kong, Kowloon Tong, Hong Kong, China}
}

\bigskip


\begin{abstract}
When a plane shock hits a two-dimensional wedge head on, it experiences a reflection-diffraction process,
and then a self-similar reflected shock moves outward as the original shock moves forward in time.
The experimental, computational, and asymptotic analysis has indicated that various patterns occur, including regular
reflection and Mach reflection.
The von Neumann conjectures on the
transition from regular to Mach reflection
involve the existence, uniqueness, and stability of
regular shock reflection-diffraction configurations,
generated by concave cornered wedges for compressible flow.
In this paper, we discuss some recent developments in the study of the von Neumann conjectures.
More specifically, we present our recent results of the uniqueness and stability of regular shock
reflection-diffraction configurations governed by the potential flow equation in an appropriate class of solutions.
We first show that the transonic shocks in the global solutions obtained
in Chen-Feldman \cite{cf-book2014shockreflection} are convex.
Then we establish the uniqueness of global shock reflection-diffraction configurations with convex transonic shocks for any wedge angle larger than the
detachment angle or the critical angle.
Moreover, the stability of the solutions with respect to the wedge angle is also shown.
Our approach also provides an alternative way of proving the existence of the
admissible solutions established first in \cite{cf-book2014shockreflection}.
\end{abstract}

\section{Introduction}\label{sec:introduction}
We survey some recent developments in the mathematical analysis of the shock reflection-diffraction problem for potential flow
and the corresponding von Neumann conjectures on the
existence, uniqueness, and stability of
regular shock reflection-diffraction configurations
for the transition from regular to Mach reflection.
The shock reflection-diffraction problem is a lateral Riemann problem
and has been not only longstanding open in fluid mechanics
but also fundamental in the mathematical theory of multidimensional conservation laws.

\smallskip
When a planar shock hits a concave cornered wedge, the incident
shock interacts with the wedge, leading to the occurrence of
shock reflection-diffraction
({\it cf.} \cite{Chapman,VanDyke}).
Beginning from the work of E. Mach \cite{Mach} in 1878,
various patterns of shock reflection-diffraction configurations have been
observed experimentally and later numerically,
including regular reflection and Mach reflection.
The existence of the regular reflection solutions for potential flow has been now
fully understood mathematically
(see \cite{cf-annofmath20101067regularreflection, cf-book2014shockreflection}),
by reducing the shock reflection-diffraction problem to a
{\em free boundary problem}, where the unknown reflected shock
is regarded as a free boundary.
Then a natural followup fundamental question is to study the uniqueness
and stability of the solutions
we have obtained.

For the uniqueness problem, it is necessary to restrict to a class of solutions.
Recent results \cite{ChiodaroliDeLellisKreml,ChiodaroliKreml,FeireislKlingenbergKremlMarkfelder,MarkfelderKlingenberg}
show the non-uniqueness of solutions with planar shocks in the class of entropy solutions with shocks
of the Cauchy problem (initial value problem) for the multidimensional compressible Euler equations
(isentropic and full).
Our setup is different -- the problem for solutions with shocks for
potential flow is on the domain with
boundaries, so these non-uniqueness results do not apply directly.
However, these results indicate that it is natural to study the uniqueness and stability problems
in a more restricted class of solutions.
In this paper, we show the uniqueness in the class of self-similar solutions
of regular shock reflection-diffraction configurations
with convex transonic shocks, which are called admissible solutions;
see the detailed definition in \S 3.
Technically, restricting to the class
of admissible solutions allows us to reduce
the uniqueness problem for shock reflection-diffraction to
a corresponding uniqueness problem for solutions of a
free boundary problem for a nonlinear elliptic equation, which is degenerate
for the supersonic case (see Fig. \ref{figure: free boundary problems-1} below).

A key property of admissible solutions which we employ in the uniqueness proof
is that the admissible solutions converge to the unique normal reflection
solution as the wedge angle tends to $\frac\pi 2$.
Then the outline of the
uniqueness argument is the following:
If there are two different admissible solutions, defined by the potential
functions $\varphi$ and $\varphi^*$,
for some wedge angle
$\theta_{\rm w}^*< \frac\pi 2$, then it suffices to:
\begin{itemize}
	\item[(i)] construct continuous families of solutions parametrized by the wedge
	angle $\theta_{\rm w} \in [\theta_{\rm w}^*, \frac\pi 2]$,
	starting from $\varphi$ and $\varphi^*$, respectively, in an appropriate norm;
	
	\smallskip
	\item[(ii)] prove {\it local uniqueness}: If two admissible
	solutions for the same wedge angle
	are close in the
    norm mentioned above, then they must be equal.
\end{itemize}
Combining this with the fact that both families converge to the unique normal
reflection as $\theta_{\rm w}\to \frac\pi 2-$, we conclude a contradiction.

Therefore, it remains to perform the two steps described above.
Both steps can be achieved if we linearize the free boundary
problem around an admissible solution, and then show that the
linearization is sufficiently regular so that
the solutions for close wedge angles can be constructed
by the implicit function theorem.
Indeed, this approach works for one regular shock reflection-diffraction
case -- the subsonic-away-from-sonic case (see \S 5 for more details).

However, it turns out that the linearization does not have such properties for
the other case -- the supersonic
case,
owing to the elliptic degeneracy near the sonic arc and relatively lower
regularity of admissible solutions
near the corner point
between the shock and the sonic arc.
For this case, instead, we develop a nonlinear approach:
We prove directly the local uniqueness property and employ it to
perturb any given admissible solution $\varphi$ for the wedge angle $\theta_{\rm w}$,
that is, to construct an admissible solution close to $\varphi$ for all wedge angles
that are sufficiently close to $\theta_{\rm w}$ by using the Leray-Schauder
degree argument in a
{\it small iteration set}.
We note that, in \cite{cf-book2014shockreflection},
the solutions have also been constructed by the Leray-Schauder degree argument, but in
a {\it large} iteration set, {\it i.e.} a subset in a space determined
by some weighted and
scaled $C^{k, \alpha}$ norms,
with bounds by the constants sufficiently large so that the {\it apriori} estimates
of the admissible solutions assure that a fixed point of the iteration map
does not occur at the boundary of the iteration set.
In the present case of
{\it small} iteration set, the similar property is shown by using the local
uniqueness.

Our proof of the local uniqueness
is based on the convexity of the reflected-diffracted transonic shock,  established in
Chen-Feldman-Xiang \cite{chenfeldmanxiangConvexity}.
We note that the convexity of the shocks is consistent with physical experiments and
numerical simulations; see e.g.
\cite{Ben-Dor,Chapman,Demyanov,Deschambault,Glaz1,Glaz2,Glaz3,Glimm,Hindman,Ivanov,Kutler,Schneyer,Serre,Shankar,Woodward},
and the references therein.
Also see \cite{CCY1,CCY2,Kurganov,Lax,LiuLax,Schulz,Serre}
for the convexity of transonic shocks in numerical Riemann solutions of the Euler equations for compressible
fluids. Mathematically, the Rankine-Hugoniot conditions on the shock whose location is unknown, together with the
nonlinear equation in the elliptic and hyperbolic regions, enforce a restriction to possible geometric
shapes of the transonic shock.
Moreover, one of our observations is that the convexity of transonic shocks is not a local property.
In fact, the uniform convexity is a result of the interaction of the cornered wedge
and the incident shock, since the reflected shock remains flat when the wedge
is a flat wall for the normal shock reflection.
In addition, for this problem, it seems to be difficult to apply the methods directly as
in \cite{Caffarelli1,Caffarelli2,Caffarelli3,Dolbeault,Plotnikov},
owing to the difference and the more complicated
structure of the boundary conditions.

In \cite{chenfeldmanxiangConvexity}, we have developed an approach in which the
global properties of solutions are incorporated in the proof of the convexity of transonic shocks.
In particular, we have introduced a general set of conditions and employed the
approach to prove the convexity of
transonic shocks under these conditions.
As a direct corollary, we have proved the
uniform convexity of transonic shocks in the two longstanding fundamental shock
problems -- the shock reflection-diffraction problem by wedges
and the
reflection problem for supersonic flows past solid ramps.

Moreover, as a byproduct of our uniqueness proof, we have
developed a new way of establishing the existence of global solutions of the shock reflection-diffraction problem
up to the detachment angle or the critical angle, based on the fine
convexity structure. Our approach is also helpful for other related mathematical
problems including free boundary problems with degeneracy.

The previous works on unsteady
flows with shocks in
self-similar coordinates include the following:
The problem of shock reflection-diffraction by a concave cornered wedges
for potential flow
has been systematically analyzed
in Chen-Feldman \cite{cf-annofmath20101067regularreflection, cf-book2014shockreflection}
and Bae-Chen-Feldman \cite{bcf-invent2009174505optimalregularity},
where the existence of regular
shock reflection-diffraction configurations has been established up to
the detachment wedge angle or the critical angle for potential flow.
For the Mach reflection, S. Chen \cite{chens-JAMSFlatMachConfigurationInPseudoStationaryFlow200863}
proved the local stability of flat Mach configuration in self-similar coordinates.
Also see \cite{ckk-maa2000313UTSDequations, ckk-SIAMjma20061947,cdx-1,DengXiang,Kim} for the unsteady
transonic small disturbance equation
and the nonlinear wave system, \cite{d-Serre} for the Chaplygin gas,
and \cite{z-actamathsinica200622177pressuregradientsystem} for the pressure-gradient system.
Meanwhile, other problems have been tackled. For the
shock diffraction problem, Chen-Feldman-Hu-Xiang \cite{CFHX} showed that
regular shock configurations cannot exist for potential flow.
For
supersonic flow past a solid ramp,
Elling-Liu \cite{ellingliu-CPAMPrandtlMeyerReflection20081347}
obtained
a first rigorous unsteady result
under certain assumptions for potential flow.
Then Bae-Chen-Feldman \cite{baechenfeldman-prandtlmayerReflection,BCF}
succeeded to remove the assumptions in \cite{ellingliu-CPAMPrandtlMeyerReflection20081347}
and established the existence theorem for global shock reflection configurations
so that the steady supersonic weak shock solution as the long-time behavior of an unsteady flow
for all physical parameters,  via new mathematical techniques
developed first in Chen-Feldman \cite{cf-book2014shockreflection}.
See also \cite{Chen,CCF,CCF2,Fang,FangXiang} and the references therein for the steady transonic
shocks over two-dimensional wedges.

The organization of this paper is the following:
In \S \ref{sec:potential flow wquation and free boundary problem},
we introduce the free boundary problem for
the shock reflection-diffraction problem.
Then the existence and regularity results
obtained in \cite{cf-book2014shockreflection} are given
in \S \ref{sec:existence}.
In \S \ref{sec:Con}, we describe the result and present the main steps in the proof of the convexity of the regular
reflected-diffracted transonic shock based on \cite{chenfeldmanxiangConvexity}.
In \S \ref{sec:local uniqueness}, we discuss our recent result and outline the proof
on the uniqueness and stability of regular shock reflection-diffraction configurations.

\section{The Potential Flow Equation and the Shock Reflection-Diffraction
	Problem}\label{sec:potential flow wquation and free boundary problem}

In this section we formulate the shock reflection-diffraction problem
as a free boundary problem for the potential flow equation in the self-similar coordinates.

\subsection{The potential flow equation}\label{subsec:potential flow rh condition}
The Euler equations
for potential flow  consist of the conservation
law of mass and Bernoulli's law:
\begin{align}
	\label{equ:1}
	&\partial_t\rho+\gr_{\mathbf{x}}\cdot(\rho\mathbf{v})=0,\\
	\label{equ:2}
	&\p_t\Phi+\frac{1}{2}|\gr_{\mathbf{x}}\Phi|^2+i(\r)=B_0,
\end{align}
where $\rho$ is the density, $\Phi$ is the velocity potential
so that $\mathbf{v}=\gr_{\mathbf{x}}\Phi$, $B_0$ is the Bernoulli constant determined by the incoming
flow and/or boundary conditions,
$\mathbf{x}=(x_1,x_2)\in\mr^2$, and
$
i(\r)=\int^{\rho}_1\frac{p'(s)}{s} {\rm d}s
$
for the pressure function $p=p(\rho)$.
For polytropic gas, by scaling,
$$
p(\r)=\frac{\r^{\g}}{\g},\quad c^2(\r)=\r^{\g-1},\quad
i(\r)=\frac{\r^{\g-1}-1}{\g-1},\quad \g>1,
$$
where $c(\r)$ is the sound speed.

The system is invariant under
the self-similar scaling:
$$
(\mathbf{x},t)\rightarrow(\alpha \mathbf{x},\alpha t),\quad
(\rho,u,v,\Phi)\rightarrow(\rho,u,v,\frac{\Phi}{\alpha})\qquad\,\,
\mbox{for $\alpha\neq0$}.
$$
Thus, we can seek self-similar solutions of the form:
$$
(\rho,u,v)(\mathbf{x},t)=(\rho,u,v)(\boldsymbol{\xi}),\quad
\Phi(\mathbf{x},t)=t\big(\varphi(\boldsymbol{\xi})+\frac{1}{2}|\boldsymbol{\xi}|^2\big)
\qquad\,\,\mbox{for $\boldsymbol{\xi}=(\xi_1,\xi_2)=\frac{\mathbf{x}}{t}$},
$$
where $\varphi$ is called a pseudo--velocity potential that satisfies
$\nabla_{\bx}\varphi=(u-\xi_1,v-\xi_2)=(U,V)$ which
is called a pseudo-velocity.
Then the pseudo--potential function
$\varphi$ satisfies the following equation for self--similar
solutions:
\begin{equation}
	\label{1.1} \mbox{div}(\rho D\varphi)+2\rho=0,
\end{equation}
where the density function $\rho=\rho(|D\varphi|^2,\varphi)$ is determined by
\begin{equation}\label{1.2}
	\rho(|D\varphi|^2,\varphi) =\big(\rho_0^{\gamma-1}-(\gamma-1)(\varphi+\frac{1}{2}|D\varphi|^2)\big)^{\frac{1}{\gamma-1}},
\end{equation}
and the divergence $\div$ and gradient $D$ are with respect to the
self--similar variables $\bx$, and $\rho_0$ is a positive constant (to be given in Problem 2.1 below) so that
$\rho_0^{\gamma-1}=(\gamma-1)B_0+1$.
Therefore, the potential function $\varphi$ is governed by the following second-order potential flow equation:
\begin{equation}\label{equ:potential flow equation}
	\div\big(\r(|D\varphi|^2,\varphi)D\varphi\big)+2\r(|D\varphi|^2,\varphi)=0.
\end{equation}
Equation \eqref{equ:potential flow equation} is a second-order
equation of mixed elliptic-hyperbolic type: It is elliptic if
and only if $|D\varphi|<c(|D\varphi|^2,\varphi)$, which is equivalent to
\begin{equation}\label{criterion:ellipticity}
	|D\varphi|<c_{\star}(\varphi,\g)\deq\sqrt{\frac{2}{\g+1}\big(\rho_0^{\gamma-1}-(\gamma-1)\varphi\big)}.
\end{equation}

If $\rho$ is a constant, then \eqref{1.1}--\eqref{1.2} imply that the corresponding
pseudo-velocity potential $\varphi$ is of the form:
$$
\varphi(\bx)=
-\frac{1}{2}|\bx|^2+ (u,v)\cdot \bx+k
$$
for constants $u$, $v$, and $k$. Such a solution is called a uniform or constant
state.

\subsection{Weak solutions and the Rankine-Hugoniot conditions}

Since shocks are involved in the problem under consideration, we define the notion
of weak solutions of equation \eqref{equ:potential flow equation}, which
admits the shocks.

\begin{definition}\label{def:weak solution}
	A function $\varphi\in W^{1,1}_{\rm loc}(\O)$ is called a weak solution of
	\eqref{equ:potential flow equation} if
	
	\bin
	\item[\rm (i)] $\rho_0^{\gamma-1}-\varphi-\frac{1}{2}|D\varphi|^2\geq0\quad\text{{\it a.e.} in }\O$,
	
	\smallskip
	\item[\rm (ii)] $(\r(|D\varphi|^2,\varphi),\r(|D\varphi|^2,\varphi)|D\varphi|)\in(L^1_{\rm loc}(\O))^2$,
	
	\smallskip
	\item[\rm (iii)] For every $\z\in C^{\infty}_{c}(\O)$,
	$$
	\int_{\O}\big(\r(|D\varphi|^2,\varphi)D\varphi\cdot
	D\z-2\r(|D\varphi|^2,\varphi)\z\big)\mathrm{d}\bx=0.
	$$
	\ein
\end{definition}

For a piecewise smooth solution $\varphi$ divided by a shock, it is easy to
verify that $\varphi$ satisfies the conditions in
Definition \ref{def:weak solution} if and only if $\varphi$ is a classic solution of
\eqref{equ:potential flow equation} in each smooth subregion and satisfies the following
Rankine-Hugoniot conditions across the shock:
\begin{eqnarray}
	&&[\r(|D\varphi|^2,\varphi)D\varphi\cdot\bn]_{S}=0, \label{RH conditon involve normal direction}\\[1.5mm]
	&&[\varphi]_{S}=0,  \label{RH conditon continous}
\end{eqnarray}
where $\bn$ is a unit normal to $S$.
Condition \eqref{RH conditon involve normal direction} is due to the conservation of mass, while
condition \eqref{RH conditon continous} is due to the irrotationality.

There are fairly many weak solutions to the given conservation
laws. The physically relevant solutions must satisfy the entropy condition.
For potential flow, the discontinuity of $D\varphi$ satisfying the Rankine-Hugoniot
conditions \eqref{RH conditon involve normal direction}--\eqref{RH conditon continous}
is called a shock if it satisfies the following
{\em entropy condition}: {\it The density $\r$ increases across a shock
	in the pseudo--flow direction}. From \eqref{RH conditon involve
	normal direction}, the entropy condition indicates that the normal
derivative function $\varphi_{\bn}=D\varphi\cdot\bn$ on a shock always decreases across
the shock in the pseudo--flow direction.

\subsection{Shock reflection-diffraction problem.}
The incident shock separates two constant states: state $(0)$
with density $\rho_0$ and velocity $\mathbf{v}_0=(0,0)$ ahead of the shock, and
state (1) with density $\rho_1$ and velocity $\mathbf{v}_1=(u_1,0)$ behind the shock,
where the entropy condition holds: $\rho_1>\rho_0$ on the shock.
The incident shock moves from the left to the right and hits the vertex of wedge:
$$
W:=\{\mathbf{x}\,:\, |x_2|< x_1 \tan\theta_{\rm w}, x_1>0\}
$$
at the initial time.
The slip boundary
condition ${\bf v}\cdot\bn=0$ is prescribed on the solid wedge boundary.

Then the shock reflection-diffraction problem can be formulated as follows:

\begin{problem}[Initial-boundary value problem]\label{PROBLEM-1}
	Seek a
	solution of system \eqref{equ:1}--\eqref{equ:2}
	for $B_0=\frac{\rho_0^{\gamma-1}-1}{\gamma-1}$
	with the initial condition at $t=0${\rm :}
	\begin{equation}\label{initial-condition}
		(\rho,\Phi)|_{t=0} =\begin{cases}
			(\rho_0, 0) \qquad\,\, &  \mbox{for}\,\, |x_2|>x_1\tan\theta_{\rm w}, \; x_1>0,\\[2mm]
			(\rho_1, u_1 x_1) \qquad\,\, &\mbox{for}\,\, x_1<0,
		\end{cases}
	\end{equation}
	and the slip boundary condition along the wedge boundary $\partial
	W${\rm :}
	\begin{equation}\label{boundary-condition}
		\nabla_{\mathbf{x}}\Phi\cdot \bn|_{\partial W \times {\mathbb{R}}_+}=0,
	\end{equation}
	where $\bn$ is the exterior unit normal to $\partial W$.
\end{problem}

The initial-boundary value problem, Problem \ref{PROBLEM-1}, is a lateral Riemann problem
with boundary $\partial W \times {\mathbb{R}}_+$ in the $(\mathbf{x}, t)$--coordinates.
Since state (1) does not satisfy the slip boundary condition, the solution must
differ from state (1) behind the shock
so that the shock reflection-diffraction configurations occur.
These configurations are self-similar, so
the problem can be reformulated
in the self-similar coordinates $\boldsymbol{\xi}=(\xi_1,\xi_2)$.
Depending on the data, there may be various patterns of shock
reflection-diffraction configurations,
including regular reflection and Mach reflection.

By the symmetry of the
problem with respect to the $\xi_1$--axis,
we consider only the upper half-plane
$\{\xi_2>0\}$ and prescribe the condition $\varphi_{\bn}=0$
on the symmetry line $\{\xi_2>0\}$.
Note that state (1) satisfies this condition.

We study self-similar solutions of Problem 2.1. Thus we give a formulation for the solution of Problem 2.1
in the self-similar coordinates $\boldsymbol{\xi}=(\xi_1, \xi_2)$. Let
$$
\Lambda=\mR^2_+\setminus\{\boldsymbol{\xi}\; : \; \xi_1>0,\, 0<\xi_2<\xi_1\tan\theta_{\rm w}\},
$$
where $\mR^2_+=\mR^2\cap\{\xi_1>0\}$.
Then, following Definition \ref{def:weak solution}, we have

\begin{definition}\label{def:weak solutionRegRefl}
	$\varphi\in C^{0,1}(\overline\Lambda)$ is a weak solution of the shock reflection-diffraction problem if $\varphi$ satisfies	
	equation \eqref{equ:potential flow equation} in $\Lambda$ in the weak sense, the boundary condition{\rm :}
	\begin{equation}\label{BCfor RegRefl}
		\partial_{\bn} \varphi=0\qquad \mbox{on }\partial\Lambda,
	\end{equation}
	and the asymptotic condition{\rm :}
	\begin{equation}\label{2.13a}
		\displaystyle
		\lim_{R\to\infty}\|\varphi-\overline{\varphi}\|_{0, \Lambda\setminus
			B_R(0)}=0,
	\end{equation}
	where
	\begin{equation*}
		\bar{\varphi}=
		\begin{cases} \varphi_0 \qquad\mbox{for}\,\,\,
			\xi_1>\xi_1^0,\,\,\, \xi_2>\xi_1 \tan\theta_{\rm w},\\[1mm]
			\varphi_1 \qquad \mbox{for}\,\,\,
			\xi_1<\xi_1^0, \,\,\, \xi_2>0,
		\end{cases}
	\end{equation*}
	and $\xi_1^0>0$ is the location of the incident shock.
\end{definition}

\subsection{Solutions of regular reflection structure}
\label{sectShockReflDiffrProblem}

We will show that, for certain values of parameters, there exist self-similar solutions of the regular reflection structure
for the shock reflection-diffraction problem and, moreover, these solutions are unique in the class of self-similar solutions
of such a structure.

\begin{figure}[htp]
\begin{center}
	\begin{minipage}{0.50\textwidth}
		\centering
		\includegraphics[width=0.65\textwidth]{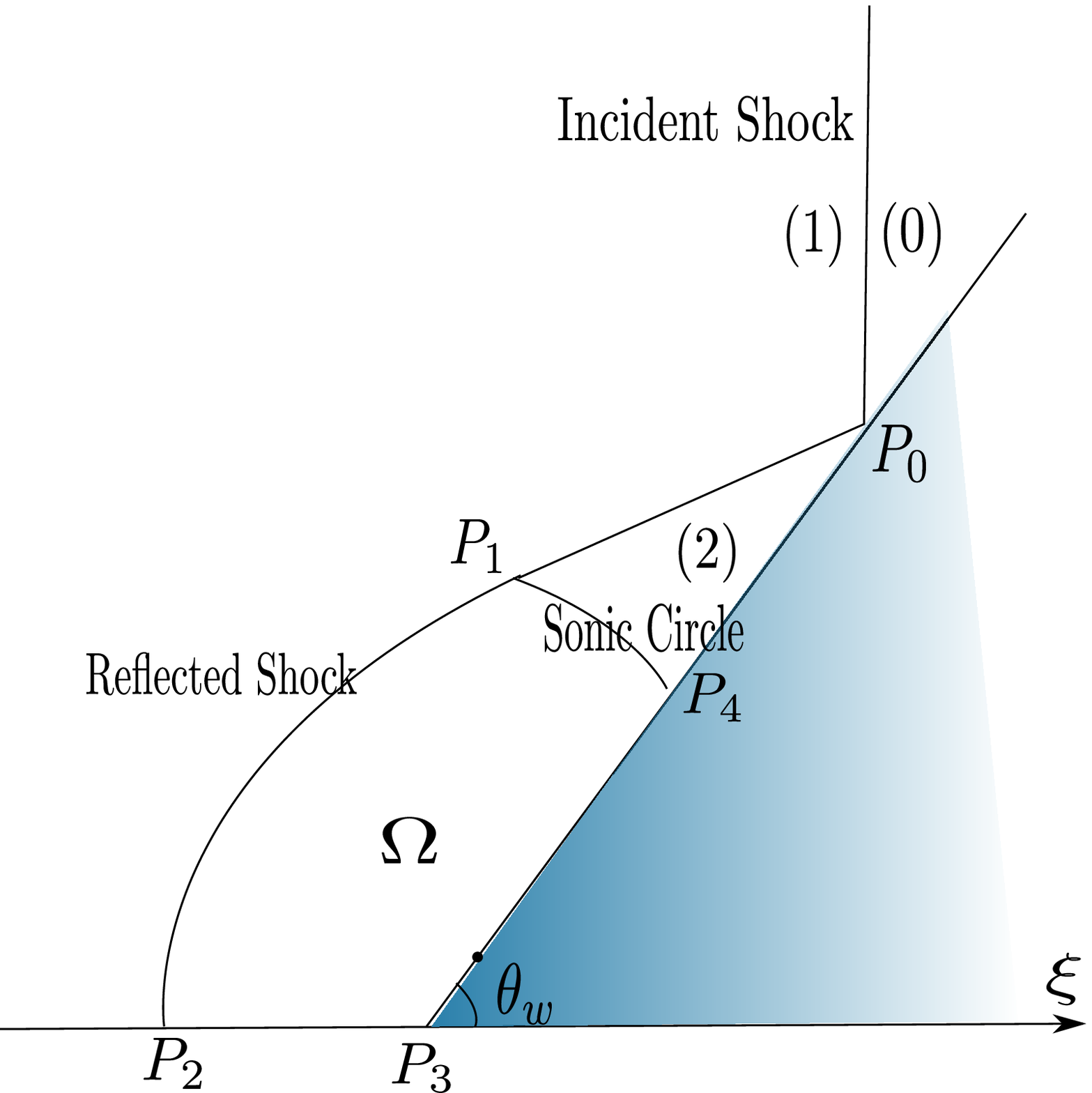}
		\caption{Supersonic regular shock reflection-diffraction configuration}
		\label{figure: free boundary problems-1}
	\end{minipage}
	\hspace{-0.5in}
	\begin{minipage}{0.50\textwidth}
		\centering
		\includegraphics[width=0.65\textwidth]{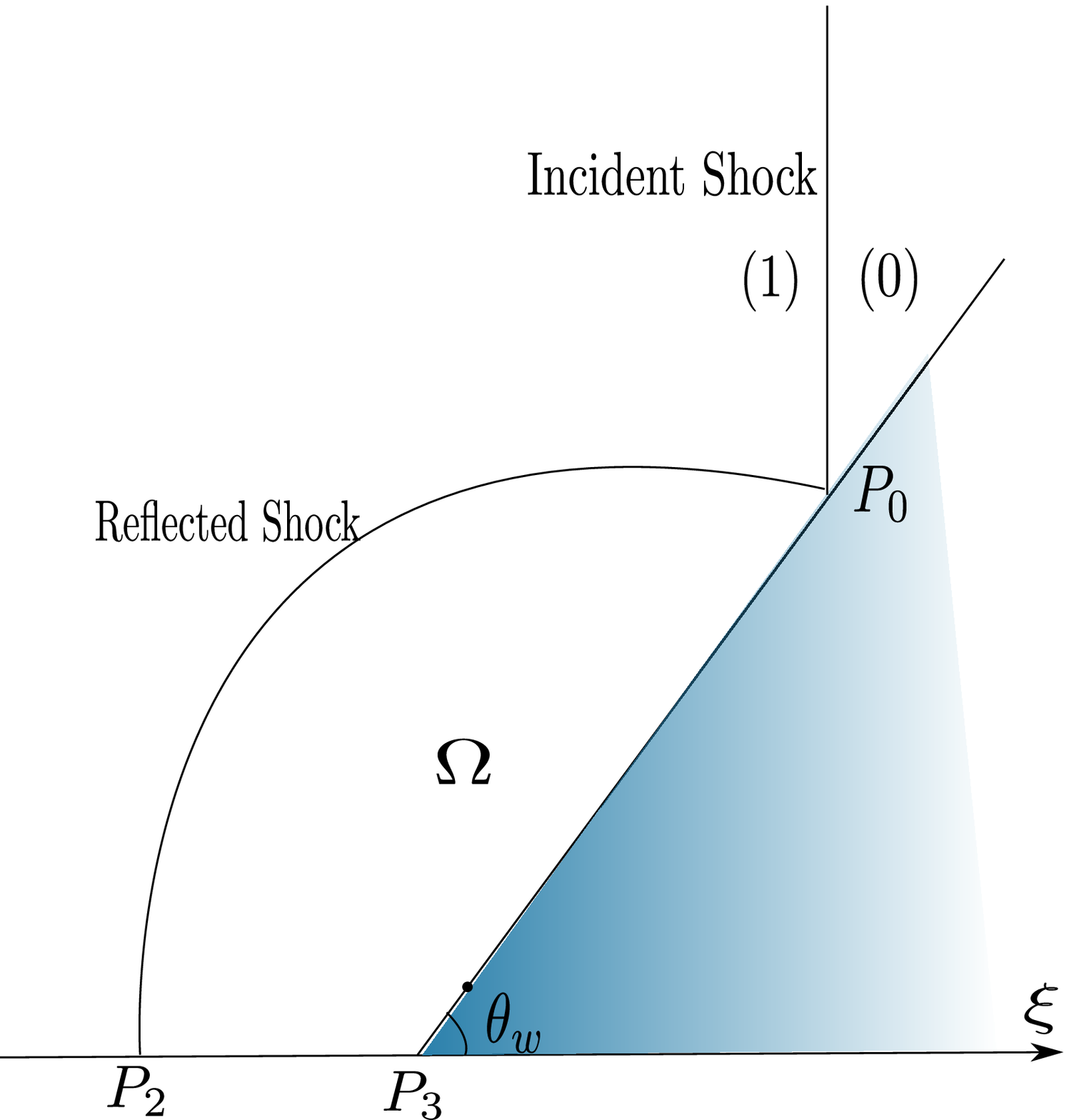}
		\caption{Subsonic$\,\,\,$  regular shock reflection-diffraction configuration}
		\label{figure: free boundary problems-2}
	\end{minipage}
\end{center}
\end{figure}

Figs. \ref{figure: free boundary problems-1}--\ref{figure: free boundary problems-2}
show two different regular shock reflection-diffraction configurations in the self-similar
coordinates.
The regular reflection solutions are piecewise smooth; more precisely,
they are smooth away from the incident and reflected-diffracted shocks,
as well as the sonic circle (which is a weak discontinuity) for the supersonic
case as shown in Fig. \ref{figure: free boundary problems-1}.

A necessary condition for the existence of piecewise-smooth regular shock reflection-diffraction
configurations is the existence of the constant state $(2)$ with the
pseudo-potential $\varphi_2$ that satisfies both
the slip boundary condition on the wedge and the Rankine-Hugoniot conditions
with state $(1)$ across the flat shock $S_1=\{\varphi_1=\varphi_2\}$, which passes
through point $P_0$ where the incident shock meets the wedge boundary.
Therefore, it requires the following three conditions at $P_0$:
\begin{equation}\label{condState2}
	\begin{split}
		&D\varphi_2\cdot\bn_{\rm w}=0,\\
		&\varphi_2=\varphi_1,\\
		&\r(|D\varphi_2|^2,\varphi_2)D\varphi_2\cdot\bn_{S_1}=\rho_1D\varphi_1\cdot\bn_{S_1},\;\;
	\end{split}
\end{equation}
where  $\bn_{\rm w}$ is the outward normal to the wedge boundary,
and $\bn_{S_1}=\frac{D(\varphi_1-\varphi_2)}{|D(\varphi_1-\varphi_2)|}$.

It is well-known (see {\it e.g.} \cite{cf-book2014shockreflection}) that,
for given parameters $(\rho_0, \rho_1)$ of states $(0)$ and $(1)$, there
exists a detachment angle $\theta_{\rm w}^{\rm d}\in (0, \frac{\pi}{2})$ such
that the algebraic equations (\ref{condState2}) have two solutions for each wedge angle
$\theta_{\rm w}\in (\theta_{\rm w}^{\rm d}, \frac{\pi}{2})$, which become
equal when $\theta_{\rm w}=\theta_{\rm w}^{\rm d}$.
Then two two-shock configurations occur at $P_0$ when
$\theta_{\rm w}\in (\theta_{\rm w}^{\rm d}, \frac\pi 2)$.
For each $\theta_{\rm w}$,
state $(2)$ with the smaller density is called a weak state $(2)$.
In this paper, state $(2)$ always refers to the weak one, since the strong state $(2)$ is ruled out
by the stability/continuity criterion as introduced first by Chen-Feldman in
\cite{cf-annofmath20101067regularreflection}; see also \cite{cf-book2014shockreflection}.
Depending on the wedge angle,  state $(2)$ can be either supersonic or
subsonic at $P_0$. Moreover, for $\theta_{\rm w}$ near $\frac\pi 2$
(resp. for $\theta_{\rm w}$ near $\theta_{\rm w}^{\rm d}$),
state $(2)$ is supersonic
(resp. subsonic) at $P_0$.
The type of state $(2)$ at $P_0$ determines the type of reflection, \emph{i.e.} supersonic or subsonic,
as shown  in Figs. \ref{figure: free boundary problems-1}--\ref{figure: free boundary problems-2}.

We consider solutions of the structure shown in Figs.
\ref{figure: free boundary problems-1}--\ref{figure: free boundary problems-2}.
Outside of region $\Omega$, the flow consists of the uniform states $(0)$, $(1)$, and $(2)$
as indicated on
the pictures, separated by the straight shocks. In particular, the incident
shock separating states (0) and (1) within $\Lambda$ is the half-line
$S_0^+=\{(\xi_1, \xi_2)\; : \; \xi_1=\xi_{1,P_0},\, \xi_2>\xi_{2, P_0}\}$.
The flow is non-uniform
and pseudo-subsonic in $\Omega$.
Here  $\Omega$ is an open bounded connected domain,
and $\partial \Omega = \overline\Gsh\cup\overline\Gso\cup\overline\Gw\cup
\overline{\Gamma_{\rm sym}}$, where curve
$\Gsh$ with endpoints $P_1$ and $P_2\in\{\xi_2=0\}$ in the supersonic case
(resp. $P_0$ and $P_2\in\{\xi_2=0\}$ in the subsonic case) is a transonic shock  which separates a constant
state $(1)$ outside $\Omega$ from a
pseudo-subsonic (non-constant) state inside $\Omega$,
and $\overline\Gso\cup\overline\Gw\cup\overline{\Gamma_{\rm sym}}$ is the fixed boundary
with  arc  $\Gso$ between points $P_1$ and $P_4$ of
the pseudo-sonic circle of state $(2)$ (we also use notation
$\overline\Gso=\{P_0\}$ for
the subsonic reflection case as shown in Fig. \ref{figure: free boundary problems-2}),
the line segment $\Gw$ is the part of  $\partial\Omega$ on the wedge boundary, {\it i.e.}
$\Gw=P_3P_4$ in the supersonic case
and $\Gw=P_0P_3$ in the subsonic case,
and $\Gamma_{\rm sym}=P_2P_3$ is
the part of $\partial\Omega$ on the symmetry line $\{\xi_2=0, \xi_1<0\}$.

\section{Existence and Regularity of Regular Shock Reflection-Diffraction
	Configurations}\label{sec:existence}

We first notice that a key obstacle to the existence of
regular shock reflection-diffraction configurations is an additional possibility
that, at the critical wedge
angle $\theta_{\rm w}^{\rm c}\in(\theta^{\rm d}_{\rm w},\frac{\pi}{2})$,
the reflected shock
$P_0P_2$ may attach to the wedge vertex $P_3$, {\it i.e.} $P_2=P_3$.
We can rule out such a solution if $u_1\le c_1$. In the opposite case
$u_1> c_1$, there would be a possibility that the reflected shock is
attached to the wedge vertex, as the experiments show ({\it e.g.} \cite[Fig. 238]{VanDyke}).
We note that the condition on $(u_1, c_1)$ can be explicitly expressed through
parameters $(\rho_0, \rho_1)$ of states (0) and (1), besides $\gamma\ge 1$,
by using \eqref{1.2} and the
Rankine-Hugoniot conditions
on the incident shock.
Recall that $\rho_1>\rho_0$.
It can be shown that there exists $\rho^c>\rho_0$ such that
$$
u_1\le c_1\;\;\mbox{ iff } \rho_1\in[\rho_0, \rho^c],\qquad\,\,
u_1> c_1\;\;\mbox{ iff }\rho_1\in[\rho^c, \infty).
$$

Now we state the existence and regularity results for the solutions of
shock reflection-diffraction problem which have the regular reflection structure as
on Fig. {\rm \ref{figure: free boundary problems-1}}--{\rm \ref{figure: free boundary problems-2}},
established in
Chen-Feldman \cite{cf-book2014shockreflection}. We prove these results
in the class of {\em admissible solutions} of the regular reflection problem, defined as follows:
\begin{definition}\label{admisSolnDef}
	Let $\theta_{\rm w}\in (\theta_{\rm w}^{\rm d},\frac{\pi}2)$.
	A function $\varphi\in C^{0,1}(\overline\Lambda)$ is an admissible solution of the regular reflection
	problem \eqref{equ:potential flow equation} and \eqref{BCfor RegRefl}--\eqref{2.13a}
	if $\varphi$ is a solution in the sense of Definition {\rm \ref{def:weak solutionRegRefl}}
	and satisfies the following properties:
	\begin{enumerate}[\rm (i)]
		\item\label{RegReflSol-Prop0}
		The structure of solutions is as follows{\rm :}
		
		\smallskip
		\begin{itemize}
			\item
			If $|D\varphi_2(P_0)|>c_2$,	then $\varphi$ is of the {\em supersonic} regular shock
			reflection-diffraction configuration described in \S \ref{sectShockReflDiffrProblem}
			and shown on Fig. {\rm \ref{figure: free boundary problems-1}}			
			and satisfies{\rm :}
			
			The reflected-diffracted shock $\Gsh$ is $C^{2}$ in its relative interior.
            Curves $\Gsh$, $\Gso$, $\Gw$, and $\Gamma_{\rm symm}$ do not
			have common points except their endpoints.

            $\varphi$ satisfies the following properties:
			\begin{align}
				& \varphi\in C^{0,1}(\Lambda)\cap C^1(\Lambda\setminus (\overline{S_0^+}\cup \overline{P_0P_1P_2})),\nonumber\\
				&\varphi\in C^{1}(\overline{\Omega})\cap C^{3}(\overline\Omega\setminus(\overline\Gso\cup\{P_2, P_3\})),\nonumber\\[2mm]
			&\varphi=\left\{\begin{array}{ll}
					\varphi_0 \qquad\mbox{for}\,\, \xi_1>\xi_1^0 \mbox{ and } \xi_2>\xi_1\tan\theta_{\rm w},\\[2mm]
					\varphi_1 \qquad\mbox{for}\,\, \xi_1<\xi_1^0
					\mbox{ and above curve} \,\, P_0P_1P_2,\\[2mm]
					\varphi_2 \qquad \mbox{in region}\,\, P_0P_1P_4.
				\end{array}\right. \label{phi-states-0-1-2-MainThm}
			\end{align}

			\smallskip
			\item If $|D\varphi_2(P_0)|\le c_2$,	then $\varphi$ is of the {\em subsonic} regular shock
			reflection-diffraction configuration described in \S \ref{sectShockReflDiffrProblem}
			and shown on Fig. {\rm \ref{figure: free boundary problems-2}}			
			and satisfies:
			
			The reflected-diffracted shock $\Gsh$ is $C^{2}$ in its relative interior. Curves $\Gsh$,  $\Gw$, and $\Gamma_{\rm symm}$ do not
			have common points except their endpoints.
			
            $\varphi$ satisfies the following properties:
			\begin{align}
				& \varphi\in C^{0,1}(\Lambda)\cap C^1(\Lambda\setminus (\overline{S_0^+}\cup\overline\Gsh)),\nonumber\\
				&\varphi\in  C^{1}(\overline{\Omega})\cap
				C^3(\overline\Omega\setminus\{P_0, P_3\}), \nonumber\\[2mm]
      		&\varphi=\left\{\begin{array}{ll}
					\varphi_0 \qquad&\mbox{for}\,\, \xi_1>\xi_1^0 \mbox{ and } \xi_2>\xi_1\tan\theta_{\rm w},\\[2mm]
					\varphi_1 \qquad&\mbox{for}\,\, \xi_1<\xi_1^0
					\mbox{ and above curve} \,\, P_0P_2,\\[2mm]
					\varphi_2(P_0) \qquad &\mbox{at $P_0$},
				\end{array}\right. \label{phi-states-0-1-2-MainThm-Subs}\\[2mm]
            &D\varphi_{|\Omega}(P_0)=D\varphi_2(P_0). \nonumber\frac{}{}
            \end{align}
		\end{itemize}
		Moreover, in both supersonic and subsonic cases,
		denote $\Gsh^{\rm ext}=\Gsh\cup \{P_0\}\cup\Gsh^-$, where $\Gsh^-$ is the reflection of $\Gsh$
		with respect to the $\xi_1$-axis. Then curve $\Gsh^{\rm ext}$ is $C^1$ in its relative interior.
		
		\smallskip		
		\item\label{RegReflSol-Prop1}
		Equation \eqref{equ:potential flow equation} is strictly elliptic in
		$\overline\Omega\setminus\,\overline{\Gso}$, {\it i.e.}
		$$
		|D\varphi|<c(|D\varphi|^2, \varphi)
		\qquad\mbox{ in }\; \overline\Omega\setminus\,\overline{\Gso},
		$$
		where, for the subsonic and sonic cases,  we use notation $\overline{\Gso}=\{P_0\}$.
		\smallskip
		\item\label{RegReflSol-Prop1-1}
		$\partial_{\bn}\varphi_1>\partial_{\bn}\varphi>0$ on $\Gsh$, where $\bn$ is the normal
		to $\Gsh$, pointing to the interior of $\Omega$.
		
		\smallskip
		\item \label{RegReflSol-Prop1-1-1}
		$\varphi_2\le\varphi\le\varphi_1$ in $\Omega$.
		
		\smallskip
		\item\label{RegReflSol-Prop2}
		Let $\mathbf{e}_{S_1}$ be the unit vector parallel to $S_1:=\{\varphi_1=\varphi_2\}$, oriented so that
		$\mathbf{e}_{S_1}\cdot D\varphi_2(P_0)>0${\rm :}
		\begin{align}\label{defEs1a}
			&\mathbf{e}_{S_1}=
			-\frac {(v_2,\, u_1-u_2)}{\sqrt{(u_1-u_2)^2+v_2^2}}.
		\end{align}
		Let $\mathbf{e}_{\xi_2}=(0,1)$. Then
		\begin{equation}\label{nonstritConeMonot}
			\partial_{\mathbf{e}_{S_1}}(\varphi_1-\varphi)\leq0, \quad\,\,
			\partial_{\xi_2}(\varphi_1-\varphi)\leq0 \qquad\;\;\mbox{ on } \Gsh.
		\end{equation}	
	\end{enumerate}
\end{definition}

\begin{remark}\label{admisEquiv1-2}
	It can be shown that Definition {\rm \ref{admisSolnDef}} is equivalent to the definition of admissible
	solutions in {\rm \cite{cf-book2014shockreflection}}{\rm ;} see Definitions {\rm 15.1.1}--{\rm 15.1.2} there.
	Thus, all the estimates and properties of admissible solutions shown in {\rm \cite{cf-book2014shockreflection}}
	hold for the admissible solutions defined above.
	In particular, the admissible solutions converge (in an appropriate sense) to the normal
	reflection solution as $\theta_{\rm w}\to \frac\pi 2{-}$.
\end{remark}

\begin{remark}\label{vectorS1-rmk}
	For the supersonic case, $\mathbf{e}_{S_1}$ defined by \eqref{defEs1a}
	has the expression{\rm :}
	$$
	\mathbf{e}_{S_1}=\frac{P_1 -P_0}{|P_1- P_0|}.
	$$
	Moreover, in the supersonic {\rm (}resp. subsonic/sonic{\rm )} case,
	$\mathbf{e}_{S_1}$ is tangential to $\Gsh$ in its upper endpoint
	$P_1$ {\rm (}resp. $P_0${\rm )}, because
	$(\varphi, D\varphi)|_{\Omega}=(\varphi_2, D\varphi_2)$ at that point,
	and its orientation at that endpoint of $\Gsh$  is towards the relative interior
	of $\Gsh$.
\end{remark}

\begin{remark}\label{eqnInOmega-rmk}
	Since the admissible solution $\varphi$ is a weak solution in the sense of Definition {\rm \ref{def:weak solutionRegRefl}}
	and is of regularity as in \eqref{RegReflSol-Prop0} of Definition {\rm \ref{admisSolnDef}},
	it satisfies \eqref{equ:potential flow equation} classically in $\Omega$ with
	the Rankine-Hugoniot conditions{\rm :}
	\begin{equation}\label{equ:boundary-RH-RegRefl}
		\varphi=\varphi_1,\quad \r(|D\varphi|^2,\varphi)D\varphi\cdot\bn=\rho_1D\varphi_1\cdot\bn
		\qquad\,\,
		\mbox{on $\Gsh$},
	\end{equation}
	and the boundary condition{\rm :}
	\begin{equation}\label{equ:boundary-ofOmega-RegRefl}
		\partial_{\bn}\varphi=0\qquad\;\mbox{ on $\Gw\cup\Gamma_{\rm sym}$}.
	\end{equation}
\end{remark}

\begin{remark}\label{nonconsta-rmk}
	The admissible solution $\varphi$ is not a constant state  in $\Omega$. Indeed, if
	$\varphi$ is a constant state in $\Omega$,
	then $\varphi=\varphi_2$ in $\Omega${\rm :}
	This follows from both \eqref{phi-states-0-1-2-MainThm} for the supersonic
	case {\rm (}since $\varphi$ is $C^{1}$ across $\Gso${\rm )} and
	property
	$(\varphi, D\varphi)=(\varphi_2, D\varphi_2)$ at $P_0$ for the subsonic case.
	However, $\varphi_2$ does not satisfy \eqref{equ:boundary-ofOmega-RegRefl} on $\Gamma_{\rm sym}$
	since $\mathbf{v}_2=(u_2, v_2)=(u_2, u_2\tan\theta_{\rm w})$ with $u_2>0$
	and $\theta_{\rm w}\in(0, \frac\pi 2)$.
\end{remark}

The following theorem shows that the admissible solution has additional regularity and monotonicity properties.

\begin{theorem}[Properties of admissible solutions]\label{Regularity-RegRefl-Th1}
$\,\,$	There exits a constant $\alpha=\alpha(\rho_0,\rho_1,\gamma)\in (0, \frac{1}{2})$
	such that any admissible solution in the sense of Definition {\rm \ref{admisSolnDef}}
	with wedge angle $\theta_{\rm w}\in (\theta_{\rm w}^{\rm d},\frac{\pi}2)$ has the following properties{\rm :}
	\begin{enumerate}[\rm (i)]
		\item\label{RegReflSol-Prop0-th}
		Additional regularity{\rm :}
		\smallskip
		\begin{itemize}
			\item
			If $|D\varphi_2(P_0)|>c_2$,	i.e. when $\varphi$ is of the {\em supersonic} regular shock
			reflection-diffraction configuration
			as in Fig. {\rm \ref{figure: free boundary problems-1}}, 		
			it satisfies
			$\varphi\in C^{1,\alpha}(\overline{\Omega})\cap C^{\infty}(\overline\Omega\setminus(\overline\Gso\cup\{P_3\}))$,
			and	$\varphi$ is $C^{1,1}$ across $\Gso$, including endpoints $P_1$ and $P_4$.
			The reflected-diffracted shock $P_0P_1P_2$ is $C^{2,\beta}$
			up to its endpoints for any $\beta\in [0, \frac{1}{2})$,
			and	$C^\infty$ except $P_1$.
			
			\smallskip
			\item If $|D\varphi_2(P_0)|\le c_2$,
			i.e. when $\varphi$ is of the {\em subsonic} regular shock reflection-diffraction configuration
			as in Fig. {\rm \ref{figure: free boundary problems-2}}, it satisfies
			$$
			\varphi\in  C^{1,\beta}(\overline{\Omega})\cap
			C^{1,\alpha}(\overline\Omega\setminus \{P_0\})
			\cap
			C^{\infty}(\overline\Omega\setminus \{P_0,  P_3\})
			$$
			for some $\beta=\beta(\rho_0,\rho_1,\gamma, \theta_w)\in (0, \alpha]$ where
			$\beta$ is non-decreasing with respect to $\theta_w$,
			and the reflected-diffracted shock $\Gsh$ is $C^{1,\beta}$
			up to its endpoints
			and $C^\infty$ except $P_0$.			
		\end{itemize}

		\smallskip
		\item\label{RegReflSol-Prop2-th}
		For each $\mathbf{e}\in Con(\mathbf{e}_{S_1}, \mathbf{e}_{\xi_2})$,
		\begin{equation}\label{strictCone-regrefl}
			\partial_{\mathbf{e}}(\varphi_1-\varphi)< 0 \qquad\mbox{in }\overline\Omega,
		\end{equation}
		where the vectors  $\mathbf{e}_{S_1}$ and  $\mathbf{e}_{\xi_2}$
		 are
		defined in  Definition {\rm \ref{admisSolnDef}}\eqref{RegReflSol-Prop2}, and
\begin{equation}\label{3.5a}
			Con(\mathbf{e}_{S_1}, \mathbf{e}_{\xi_2}):=\{a\mathbf{e}_{S_1}+b\mathbf{e}_{\xi_2}\; : \;a,b>0\}.
		\end{equation}

		\smallskip
\item
		\label{RegReflSol-Prop3}
Denote by $\bn_{\rm w}$ the
unit interior normal on $\Gw$ {\rm (}with respect to $\Omega${\rm )}, {\it i.e.}
$\bn_{\rm w}=(-\sin\theta_{\rm w}, \cos\theta_{\rm w})$.
Then
$\partial_{\bn_{\rm w}}(\varphi-\varphi_2)\leq0$ in $\overline\Omega$.
\end{enumerate}
\end{theorem}

\begin{remark}
	$Con(\mathbf{e}_{S_1}, \mathbf{e}_\eta)=\{a\mathbf{e}_{S_1}+
	b\mathbf{e}_\eta:\,a,b>0\}$ is an open set{\rm ;} that is, it does not include
	the directions of
	$\mathbf{e}_{S_1}$ and $\mathbf{e}_{\xi_2}$.
\end{remark}

\smallskip

Now we state the results on the existence of admissible solutions.

\begin{theorem}[Global solutions up to the detachment angle for the case: $u_1\le c_1$]\label{ExistRegRefl-Th1}
	Let the initial data $(\rho_0, \rho_1, \gamma)$ satisfy that $u_1\le c_1$.
	Then, for each $\theta_{\rm w}\in (\theta_{\rm w}^{\rm d},\frac{\pi}2)$, there
	exists an admissible solution of the regular reflection problem
	in the sense of Definition {\rm \ref{admisSolnDef}}. Note that these solutions
	satisfy the properties stated in Theorem {\rm \ref{Regularity-RegRefl-Th1}}.
\end{theorem}

\smallskip
\begin{theorem}[Global solutions up to the detachment angle for the case: $u_1> c_1$]\label{ExistRegRefl-Th2}
	Let the initial data  $(\rho_0, \rho_1, \gamma)$ satisfy that $u_1>c_1$.	
	Then there is
	$\theta_{\rm w}^{\rm c}\in[\theta^{\rm d}_{\rm w},\frac{\pi}{2})$ such
	that, for each $\theta_{\rm w}\in (\theta_{\rm w}^{\rm c},\frac{\pi}2)$,
	there
	exists an admissible solution of the regular reflection problem
	in the sense of Definition {\rm \ref{admisSolnDef}}.
	Note that these solutions
	satisfy the properties stated in Theorem {\rm \ref{Regularity-RegRefl-Th1}}.
	
	If $\theta_{\rm w}^{\rm c}>\theta^{\rm d}_{\rm w}$, then, for the wedge angle $\theta_{\rm w}=\theta_{\rm w}^{\rm c}$, there
	exists an attached shock solution $\varphi$ with all the properties listed in
	Definition {\rm \ref{admisSolnDef}}
	and Theorem {\rm \ref{Regularity-RegRefl-Th1}(\ref{RegReflSol-Prop2-th})--(\ref{RegReflSol-Prop3})} except that $P_2=P_3$.
	In addition, for the regularity of solution $\varphi$, we have
	\begin{itemize}
		\item For the supersonic case with $\theta_{\rm w}=\theta^{\rm c}_{\rm w}$,
		$$
		\varphi\in C^\infty(\overline{\Omega}\backslash(\overline\Gso\cup\{P_3\}))\cap C^{1,1}(\overline{\Omega}\backslash\{P_3\})\cap C^{0,1}(\overline{\Omega}),
		$$
		and the reflected shock $P_1P_2P_3$ is Lipschitz up to the endpoints, $C^{2,\beta}$ with any $\beta\in[0,\frac{1}{2})$ except point $P_3$,
		and $C^{\infty}$ except points $P_1$ and $P_3$.
		
		\smallskip
		\item For the subsonic case with $\theta_{\rm w}=\theta^{\rm c}_{\rm w}$,
		$$
		\varphi\in C^\infty(\overline{\Omega}\backslash\{P_1,\,P_3\})
		\cap C^{1,\beta}(\overline{\Omega}\backslash\{P_3\})\cap C^{0,1}(\overline{\Omega})
		$$
		for $\beta$ as in Theorem  {\rm \ref{Regularity-RegRefl-Th1}},
		and the reflected shock $P_1P_2P_3$ is Lipschitz up to the endpoints, $C^{1,\beta}$ except
		point $P_3$, and $C^{\infty}$ except points $P_1$ and $P_3$.
	\end{itemize}
\end{theorem}
In the next two sections, \S 4--\S5, we show how the convexity of the transonic shocks and the uniqueness
of the admissible solutions
can be achieved.

\section{Convexity of Transonic Shocks in the Shock Reflection-Diffraction Configurations}\label{sec:Con}

We first note that, for an admissible solution,  $\Gsh$ is a graph in any direction
$\mathbf{e}\in Con:= {Con(\mathbf{e}_{S_1},\mathbf{e}_{\eta})}$,
where $Con(\mathbf{e}_{S_1},\mathbf{e}_{\eta})$ is defined in \eqref{3.5a}.
For the subsonic/sonic reflections case, we denote $P_1:=P_0$ so that $\Gsh$ has endpoints
$P_1$ and $P_2$ in all cases.
More precisely,
the following was shown in \cite{cf-book2014shockreflection}, as a consequence
of Theorem \ref{Regularity-RegRefl-Th1}(\ref{RegReflSol-Prop2-th}):

\begin{lemma}\label{shockIsGraphLemma}
Let $\varphi$ be an admissible solution. Denote $\phi:=\varphi-\varphi_1$.
Let $\bt_{P_1}$ be a unit tangent vector to $\Gsh$ at $P_1$,
	directed into the interior of $\Gsh$.
	Let $\mathbf{e}\in Con$, and let $\mathbf{e}^{\perp}$ be the orthogonal unit vector
	to $\mathbf{e}$ with $\mathbf{e}^{\perp}\cdot\bt_{P_1}>0$.
	Let $(S,T)$ be the coordinates with respect to basis $\{\mathbf{e}, \mathbf{e}^{\perp}\}$
    so that $T_{P_2}>T_{P_1}$.
	Then there exists $f_\mathbf{e}\in C^{1}(\mr)$ such that
	
	\smallskip
	\begin{enumerate}[\rm (a)]
		\item \label{lem:shock graph-i1}
		$\Gsh=\{S=f_\mathbf{e}(T)\,:\,T_{P_1}<T<T_{P_2}\}$, $\Omega\subset\{S<f_\mathbf{e}(T) :\,T\in\mr\}$,
		$P_1=(f_\mathbf{e}(T_{P_1}),T_{P_1})$, $P_2=(f_\mathbf{e}(T_{P_2}),T_{P_2})$,
		and $f_\mathbf{e}\in C^\infty(T_{P_1},T_{P_2})${\rm ;}
		
		\smallskip
		\item\label{lem:shock graph-i2}
		The directions of the tangent lines to $\Gsh$ lie between $\bt_{P_1}$ and $\bt_{P_2}${\rm ;}
		that is, in the $(S,T)$--coordinates,
		$$
		\qquad -\infty<\frac{\bt_{P_2}\cdot \mathbf{e}}{\bt_{P_2}\cdot \mathbf{e}^{\perp}}=f_\mathbf{e}'(T_{P_2})\leq
		f'_\mathbf{e}(T)\leq f'_\mathbf{e}(T_{P_1})
		=\frac{\bt_{P_1}\cdot \mathbf{e}}{\bt_{P_1}\cdot \mathbf{e}^{\perp}}<\infty
		$$
        for any $T\in(T_{P_1},T_{P_2})${\rm ;}		

		\item\label{lem:shock graph-i2-00}
		$\bn(P)\cdot \mathbf{e}<0$ for any $P\in\Gsh${\rm ;}
		
		\smallskip
		\item\label{lem:shock graph-i2-01}
		$\phi_\mathbf{e}>0$ on $\Gsh${\rm ;}
		
		\smallskip
		\item\label{lem:shock graph-i3}
		For any $T\in(T_{P_1},T_{P_2})$,
		$$
		\phi_{\bt\bt}(f_\mathbf{e}(T),T)<0\quad \Longleftrightarrow
		\quad f''_\mathbf{e}(T)>0,
		$$
		and
		$$
		\phi_{\bt\bt}(f_\mathbf{e}(T),T)>0\quad  \Longleftrightarrow
		\quad f''_\mathbf{e}(T)<0.
		$$
	\end{enumerate}
\end{lemma}

\smallskip
In \cite{chenfeldmanxiangConvexity},
we provide a framework for the convexity of transonic shocks in the self-similar coordinates.
Specifically, for the transonic shocks in the shock reflection-diffraction
configurations, we have the following theorem.

\begin{theorem}[Convexity of transonic shocks]\label{thm:con}
	If a solution of the shock reflection-diffraction problem is admissible in the sense of
	Definition {\rm \ref{admisSolnDef}},
	then its shock curve $\Gsh$ is strictly convex in the following sense{\rm :}
	For any $\mathbf{e}\in Con$, $f_\mathbf{e}$ from Lemma {\rm \ref{shockIsGraphLemma}}
	is concave on $(T_{P_1}, T_{P_2})$, and $f''_\mathbf{e}(T)<0$ for all $T\in (T_{P_1}, T_{P_2})$.
	That is, $\Gsh$ is uniformly convex on closed subsets of its relative interior.
	Moreover, for a regular reflection solution  in the sense of
	Definition {\rm \ref{def:weak solutionRegRefl}} with pseudo-potential
	$\varphi\in C^{0,1}(\Lambda)$ satisfying
    Definition {\rm \ref{admisSolnDef}}{\rm (\ref{RegReflSol-Prop0})}--\eqref{RegReflSol-Prop1-1-1},
	the shock is strictly convex if and only if
	Definition {\rm \ref{admisSolnDef}}\eqref{RegReflSol-Prop2} holds.
\end{theorem}

Now we discuss the techniques developed in \cite{chenfeldmanxiangConvexity}
by giving the main steps in the proof of Theorem \ref{thm:con}.
While the argument in \cite{chenfeldmanxiangConvexity} is for a general domain $\Omega$,
we focus here on the regular shock reflection-diffraction configurations,
in which both the solution domain $\Omega$ and the solution structure are somewhat simpler.

\bigskip
\noindent
{\it Outline of the Proof of Theorem} \ref{thm:con}:  The proof consists of eight steps,
while the first three steps are general properties of
shock reflection-diffraction solutions; see \cite{cf-book2014shockreflection}.
Below we use notation $\phi:=\varphi-\varphi_1$.

\medskip
1. We establish a relation between the extrema of the solution and the geometric shape
of the transonic shock.
For a fixed unit vector $\mathbf{e}\in\mr^2$, denote $w:=\partial_\mathbf{e}\phi$ in $\Omega$.
We show that, if a local minimum (resp. maximum) of $w$ is attained at $P\in\Gsh^{0}$
and $\bn(P)\cdot \mathbf{e}<0$,
then $\phi_{\bt\bt}>0$ (resp. $\phi_{\bt\bt}<0$) at $P$,
where $\bn$ denotes the interior unit normal on $\Gsh$ towards $\Omega$.

\smallskip
2. We establish a nonlocal relation between the values of $\phi_\mathbf{e}$
and the positions where these values are taken.
Let $\phi$ be a solution as in Theorem {\rm \ref{thm:con}}, and let $\mathbf{e}\in Con$.
We use the coordinates from Lemma \ref{shockIsGraphLemma}.
Assume that, for two different points $P=(T,f_\mathbf{e}(T))$
and $P_1=(T_1,f_\mathbf{e}(T_1))$
on $\Gsh$,
$$
f_\mathbf{e}(T)>f_\mathbf{e}(T_1)+f'_\mathbf{e}(T)(T-T_1),
\qquad f'_\mathbf{e}(T)=f'_\mathbf{e}(T_1).
$$
Then

\smallskip
\begin{enumerate}[\rm (i)]
	\item \label{distToSonicCenter-lemma-i1}
	$d(P):=\mbox{\rm dist}(O_0,L_{P})>\mbox{\rm dist}(O_0,L_{P_1})=:d(P_1)$,
	where $O_0$ is the center of sonic circle of state $(0)$, and $L_{P}$ and $L_{P_1}$ are
	the tangent lines of $\Gsh$ at $P$ and $P_1$, respectively.
	
	\smallskip
	\item  \label{distToSonicCenter-lemma-i2}
	If the unit vector $\mathbf{e}\in Con$, then
	$$
	\phi_\mathbf{e}(P)>\phi_\mathbf{e}(P_1).
	$$
\end{enumerate}

\medskip
3. We show that the shock graph is real analytic.

\smallskip
4. We now develop a minimal/maximal chain argument.
Let $\phi$ be an admissible solution, and let $\mathbf{e}\in\mR^2$.
Note that $\phi_\mathbf{e}$  satisfies the strong maximum principle in $\Omega$.
Then we can introduce the minimal (or maximal) chain as follows:

Let $E_1, E_2\in \partial\Omega$.
We say that points $E_1$ and $E_2$ are connected by a minimal {\rm (}resp.  maximal{\rm )}
chain with radius $r$ if and only if
there exist $r>0$, integer $k_1\ge 1$, and a chain of balls $\{B_{r}(C^i)\}_{i=0}^{k_1}$ such that

\smallskip
\begin{enumerate}[\rm (i)]
	\item \label{def:minimal chain-Ia}
	$C^0=E_1$, $C^{k_1}=E_2$, and $C^i\in\overline\Omega$  for $i=0,\dots, k_1${\rm ;}
	
	\smallskip
	\item \label{def:minimal chain-Ib}
	$C^{i+1}\in \overline {B_{r}(C^i)\cap\Omega}$ for $i=0,\dots, k_1-1${\rm ;}
	
	\smallskip
	\item \label{def:minimal chain-Ic}
	$\phi_\mathbf{e}(C^{i+1})=\min_{\overline{B_{r}(C^i)\cap\Omega}}\phi_\mathbf{e}<\phi_\mathbf{e}(C^i)$ \,
	{\rm (}resp. $\phi_\mathbf{e}(C^{i+1})=\max_{\overline{B_{r}(C^i)\cap\Omega}}\phi_\mathbf{e}>\phi_\mathbf{e}(C^i)${\rm )}
	\, for $i=0,\dots, k_1-1${\rm ;}
	
	\smallskip
	\item \label{def:minimal chain-Id}
	$\phi_\mathbf{e}(C^{k_1})=\min_{\overline{B_{r}(C^{k_1})\cap\Omega}}\phi_\mathbf{e}$
	\, {\rm (}resp. $\phi_\mathbf{e}(C^{k_1})=\max_{\overline{B_{r}(C^{k_1})\cap\Omega}}\phi_\mathbf{e}${\rm )}.
\end{enumerate}

\smallskip
For such a chain $\{C^i\}_{i=0}^{k_1}$,
we also use the following terminology{\rm :} The chain starts at $E_1$ and ends at $E_2$, or
the chain is from $E_1$ to $E_2$.

\smallskip
This definition does not rule out the possibility that $B_r(C^i)\cap\partial\Omega\ne\emptyset$,
or even $C^i\in\partial\Omega$, for some or all $i=0,\dots, k_1-1$.    	
The radius $r$ is a parameter in the definition of minimal or maximal chains.	
We do not fix $r$ at this point. 	
In fact,	
the radii are determined for various chains, respectively.

Then we prove the following results:

\smallskip
\begin{enumerate}[\rm (a)]
	\item\label{lem:chainsConnectedSets}
	{\it The chains with sufficiently small radius are connected sets}.
	More precisely, there exists  $r^*$ depending only on $(\rho_0, \rho_1, \gamma)$ such that,
	for any minimal or maximal chain  $\{C^i\}_{i=0}^{k_1}$ with $r\in(0, r^*]$,
	$\displaystyle\cup_{i=0}^{k_1} B_{r}(C^i)\cap\Omega$ is connected.
	
	\smallskip
	\item\label{a}
	{\it The existence of the minimal or maximal chain of  radius $r<r^*$}.
	More precisely, if $E_1\in\partial\Omega$, and $E_1$ is not a local minimum point {\rm (}resp. maximum point{\rm )}
	of $\phi_\mathbf{e}$ with respect to $\overline\Omega$, then, for any $r\in (0, r^*)$,
	there exists a minimal {\rm (}resp. maximal{\rm )} chain $\{G^i\}_{i=0}^{k_1}$ for $\phi_\mathbf{e}$ of radius $r$,
	starting at $E_1$, {\it i.e.} $G^0=E_1$.
	Moreover, $G^{k_1}\in\partial\Omega$ is a local minimum {\rm (}resp. maximum{\rm )} point of
	$\phi_\mathbf{e}$ with respect
	to $\overline\Omega$, and $\phi_\mathbf{e}(G^{k_1})<\phi_\mathbf{e}(E_1)$ {\rm (}resp. $\phi_\mathbf{e}(G^{k_1})>\phi_\mathbf{e}(E_1)${\rm )}.
	
	\smallskip
	\item\label{lem:chainsDoNotIntersect-MinMax}
	{\it The minimal and maximal chains do not intersect}.
	Specifically,  for any $\delta>0$,	there exists $r_1^*\in(0, r^*]$ such that the following holds{\rm :}	
	Let ${\mathcal C}\subset\partial\Omega$ be connected, let $E_1$ and $E_2$ be the endpoints	of ${\mathcal C}$,
	and let there be a minimal chain $\{E^i\}_{i=0}^{k_1}$ of radius $r_1\in (0, r_1^*]$, which starts at $E_1$ and ends at $E_2$.	
	If there exists $H_1\in {\mathcal C}^0={\mathcal C}\setminus \{E_1,E_2\}$ such that
	$$
	\phi_\mathbf{e}(H_1)\ge \phi_\mathbf{e}(E_1)+\delta,
	$$
	then, for any $r_2\in (0, r_1]$, any maximal chain $\{H^j\}_{j=0}^{k_2}$ of radius $r_2$	
	starting from $H_1$ satisfies $H^{k_2}\in  {\mathcal C}^0$,	
	where ${\mathcal C}^0$ denotes the relative interior of curve  ${\mathcal C}$ as before.
	
	Note that, if $H_1$ is not a local maximum point of $\phi_\mathbf{e}$ with respect
	to $\overline\Omega$, then the existence of  the maximal chain
	$\{H^j\}_{j=0}^{k_2}$ of radius $r_2$
	starting from $H_1$
	follows from result \eqref{a}.
	
	\smallskip
	\item\label{lem:chainsDoNotIntersect-MaxMin}
	{\it Result \eqref{lem:chainsDoNotIntersect-MinMax} also holds if the roles of minimal and maximal chains are interchanged}.
	For any $\delta>0$, there exists $r_1^*\in(0, r^*]$ such that the following holds{\rm :}	
	Let ${\mathcal C}\subset\partial\Omega$ be connected, and let $E_1$ and $E_2$ be the endpoints of ${\mathcal C}$.	
	Assume that there exists a maximal chain $\{E^i\}_{i=0}^{k_1}$ of radius $r_1\in (0, r_1^*]$,	
	which starts at $E_1$ and ends at $E_2$.	If there exists $H_1\in {\mathcal C}^0$ such that
	$$
	\phi_\mathbf{e}(H_1)\le \phi_\mathbf{e}(E_1)-\delta,
	$$
	then, for any $r_2\in (0, r_1]$, any minimal chain $\{H^j\}_{j=0}^{k_2}$ of radius $r_2$,	
	starting from $H_1$, satisfies that $H^{k_2}\in  {\mathcal C}^0$.
	
	\smallskip
	\item \label{lem:chainsDoNotIntersect-MinMin}
	{\it Two minimal chains do not intersect}:
	For any $r_1\in (0, r^*]$, there exists $r_2^*=r_2^*(r_1)\in (0, r^*]$
	such that the following holds{\rm :}
	Let ${\mathcal C}\subset\partial\Omega$
	be connected, and let $E_1$ and $E_2$ be the endpoints of ${\mathcal C}$.
	Assume that there
	exists a minimal chain $\{E^i\}_{i=0}^{k_1}$ of radius $r_1\in (0, r^*]$,
	which starts at $E_1$ and ends at $E_2$.
	If there exists $H_1\in {\mathcal C}^0$ such that
	$$
	\phi_\mathbf{e}(H_1)< \phi_\mathbf{e}(E_2),
	$$
	then, for any $r_2\in(0, r_2^*]$,
	any minimal chain $\{H^j\}_{j=0}^{k_2}$ of radius $r_2$,
	starting from $H_1$,
	satisfies that $H^{k_2}\in  {\mathcal C}^0$.
\end{enumerate}

\smallskip
5. Denote by $\bn_{\rm w}$ the
unit normal on $\Gw$ pointing into $\Omega$.
By \cite[Lemma 8.2.11]{cf-book2014shockreflection},
$\bn_{\rm w}\in Con(\mathbf{e}_{S_1}, \mathbf{e}_{\xi_2})$ for
any wedge angle $\theta_{\rm w}\in (\theta_{\rm w}^{\rm d},\frac{\pi}2)$.
We use   $\mathbf{e}=\bn_{\rm w}$
for the following four steps below.
We work in the corresponding $(S,T)$--coordinates
defined in Lemma \ref{shockIsGraphLemma},
 so
it suffices to prove that the graph is concave:
$$
f''_\mathbf{e}(T)\le 0 \qquad\; \mbox{for all }T\in(T_{P_1},T_{P_2}).
$$

If there exists $\hat P\in \Gsh^0$ with $f_\mathbf{e}''(T_{\hat P})>0$,	
we prove the existence of a point $C\in \Gsh^0$ such that $f_\mathbf{e}''(T_C)\ge 0$,
and $C$ is a  point of strict local minimum
of $\phi_\mathbf{e}$ along $\Gsh$ but is {\em not} a local minimum point
of $\phi_\mathbf{e}$ relative to $\overline\Omega$.

\medskip
6. Then we prove that there exists $C_1\in \Gsh^0$ such that
there is a  minimal chain with radius $r_1$ from $C$ to $C_1$.

\medskip
7. We show that the existence of points $C$ and $C_1$ described above yields a contradiction
(which implies that there is no $\hat P\in \Gsh^0$ with $f_\mathbf{e}''(T_{\hat P})>0$).
This is proved by showing the following facts:

\smallskip
\begin{itemize}
	\item Let $A_2$ be a maximum point of $\phi_\mathbf{e}$ along $\Gsh$ lying between points $C$ and $C_1$.
	Then $A_2$ is a local maximum point of $\phi_\mathbf{e}$ relative to $\Omega$,
	and there is no point between $C$ and $C_1$ such that
	the tangent line at this point is parallel to the one at $A_2$.
	
	\smallskip
	\item Between $C$ and $A_2$, or between $C_1$ and $A_2$, there exists a local minimum point $C_2$
	of $\phi_\mathbf{e}$ along $\Gsh$ such that $C_2\ne C_1$, or $C_2\ne C$, and $C_2$ is not a local minimum point
	of $\phi_\mathbf{e}$ relative to  domain $\overline\Omega$.
	
	\smallskip
	\item Then, following an argument similar to the one used above and going through several steps,
    we arrive at the situation that the endpoint of the minimal chain cannot lie anywhere on $\partial\Omega$,
    which is a contradiction.
\end{itemize}
This indicates that $f''_\mathbf{e}\le 0$ on $\Gsh$; that is, $\Gsh$ is convex.
In the rest of the argument, we prove that  $f''_\mathbf{e}< 0$ on $\Gsh^0$.

\medskip
8. Using the fact that the shock graph is real analytic, we show that, for every $P\in \Gsh^0$,	
either $f_\mathbf{e}''(T_P)<0$ or there exists an even integer $k>2$ such that $f_\mathbf{e}^{(i)}(T_P)=0$ for all $i=2,\dots,k-1$,
and $f_\mathbf{e}^{(k)}(T_P)<0$.
This shows the strict convexity of the shock, which implies that the shock does not
contain any straight segment.
The above property is equivalent to the facts that
$\partial_{\bt}^{i}\phi(P)=0$
for all $i=2,\dots,k-1$,
and $\partial_{\bt}^k\phi(P)>0$.

\medskip
9. We show the uniform convexity of $\Gsh^0$ in the sense that
$f''_\mathbf{e}(T_P)<0$
for every $P\in \Gsh^0$,  or equivalently,
$f''_\mathbf{e}(T)<0$ on $(T_{P_1}, T_{P_2})$, for some (and thus any) $\mathbf{e}\in Con$.
In fact, if it is not true, \emph{i.e.} if $\phi_{\bt\bt}=0$ at some $P_{\rm d}$, then we can obtain a contradiction by proving that
there exists a unit vector $\mathbf{e}\in\mathbb{R}^2$ such that $P_{\rm d}$ is a
local minimum point of $\phi_\mathbf{e}$ along $\Gamma^{0}_{\rm shock}$,
but $P_{\rm d}$ is not a local minimum point of $\phi_\mathbf{e}$ in $\Omega$.
Then we can construct a minimal chain for $\phi_\mathbf{e}$ connecting $P_{\rm d}$ to $C^{k_1}\in\partial\Omega$.
We show that
\begin{itemize}
	\item $C^{k_1}\notin\Gso$,
	
	\smallskip
	\item $C^{k_1}\notin\Gw\cup\Gsy$,
	
	\smallskip
	\item $C^{k_1}\notin\Gsh$.
\end{itemize}
This implies that  $\phi_{\bt\bt}>0$ on $\Gsh^0$ so that $f''_\mathbf{e}(T)<0$ on $(T_{P_1}, T_{P_2})$;
see Lemma \ref{shockIsGraphLemma}.

\section{Uniqueness and Stability of Regular Shock Reflection-Diffraction Configurations}\label{sec:local uniqueness}
\newcommand{\dist}{\mbox{dist}}
\newcommand{\IterMap}{\mathcal{F}}
\newcommand{\IterSet}{\mathcal{K}_{\varepsilon,\delta}^{(\hat u,\hat\theta_{\rm w})}}
\newcommand{\ItRg}{Q^{\rm iter}}
\newcommand{\Symm}{\Gamma_{\rm sym}}
\newcommand{\difPet}{\phi}

In this section, we discuss the uniqueness and stability of global regular shock reflection-diffraction configurations.
More specifically, we describe the results in Chen-Feldman-Xiang \cite{Unique}.

As indicated earlier, recent results \cite{ChiodaroliDeLellisKreml,ChiodaroliKreml,MarkfelderKlingenberg,FeireislKlingenbergKremlMarkfelder}
have shown the non-uniqueness of solutions with planar shocks in the class of entropy solutions with
shocks of the Cauchy problem for the multidimensional compressible Euler system.
Moreover,
the uniqueness problem for general self-similar solutions of the Euler system
is still open ({\it cf.} \cite{ChiodaroliDeLellisKreml}).
While these results  do not apply directly to our case, they indicate that
it be natural to
study the uniqueness of solutions in some more restrictive class, instead
of general time-dependent solutions ({\it i.e.} solutions of Problem 2.1), or even general
self-similar solutions as in Definition \ref{def:weak solutionRegRefl}.

In \cite{Unique}, we have established the uniqueness of regular reflection solutions
for each wedge angle in
the class of {\em admissible solutions} introduced in Definition \ref{admisSolnDef}.

\begin{theorem}[Uniqueness]\label{thm:main theorem2}
	For any wedge angle $\theta_{\rm w}\in(\theta_{\rm w}^{\rm d},\frac{\pi}{2})$ when $u_1\leq c_1$ and
	$\theta_{\rm w}\in(\theta_{\rm w}^{\rm c},\frac{\pi}{2})$ when $u_1>c_1$, any solution, satisfying
all properties
	{\rm (\ref{RegReflSol-Prop0})}--\eqref{RegReflSol-Prop1-1-1} in Definition {\rm \ref{admisSolnDef}}
	and one of the following properties{\rm :}
	\begin{itemize}
		\item[\rm (a)] the transonic shock $\Gsh$ is convex, {i.e.} domain $\Omega$ is a convex set,
		
		\smallskip
		\item[\rm (b)] condition \eqref{nonstritConeMonot} holds,
	\end{itemize}
	is unique in the class of admissible solutions.
	Moreover, such solutions are continuous with respect to the wedge angle $\theta_{\rm w}$ in the $C^1$--norm
	$($more precisely, the continuity with respect to the norm described
	in Remark {\rm \ref{distBetweenRegReflSol}} below$)$.
\end{theorem}

\begin{remark}\label{distBetweenRegReflSol}
	For an admissible solution $\varphi$ with a wedge angle $\theta_{\rm w}$,
	we define its norm  based on its restriction to $\Omega$.
	Since region $\Omega$ depends on the solution,
	we map a unit square $Q^{\rm iter}=(0,1)^2$ to $\Omega$  and
	use this mapping to define a function $u$ on $Q^{\rm iter}$, which corresponds
	to $\varphi_{|\Omega}$.
	Furthermore, the sides of square $Q^{\rm iter}$ are mapped to the boundary parts
	$\Gso$, $\Gw$, $\Gamma_{\rm sym}$, and $\Gsh$.
	The mapping depends on $(\varphi, \theta_{\rm w})$ and is invertible{\rm ;}
	that is, given a function $u$ on $Q^{\rm iter}$
	and $\theta_{\rm w}$, we can recover $\varphi$ and $\Omega$.
	Moreover, this mapping and its inverse have
	appropriate continuity properties.
	See {\rm \cite[\S 12.2 and \S 17.2]{cf-book2014shockreflection}} for the details.
	Then we define function spaces for admissible
	solutions and ``approximate admissible solutions'' in terms of the function spaces
	for the corresponding functions $u$ on  $Q^{\rm iter}$.
	The convergence of admissible solutions $\varphi^{(i)} \to \varphi^{(\infty)}$ in
	the $C^1$--norm
	as the corresponding wedge angles
	$\theta_{\rm w}^{(i)}\to \theta_{\rm w}^{(\infty)}$,
	defined in terms of convergence in an
	appropriate norm for the functions on $\ItRg$, implies
	\begin{equation}\label{distBetweenRegReflSol-Eqn}
		\begin{split}
			&\|\varphi^{(i)}\|_{C^1(\Omega^{(i)})}\le C \qquad\mbox{ for all } i, \\
			&\|\varphi^{(i)}-\varphi^{(\infty)}\|_{C^1(\overline{\Omega^{(i)}\cap \Omega^{(\infty)}})}
			+d_{\rm H}(\overline{\Omega^{(i)}}, \overline{\Omega^{(\infty)}}) \to 0
			\qquad \mbox{as $\theta_{\rm w}^{(i)}\to \theta_{\rm w}^{(\infty)}$},
		\end{split}
	\end{equation}
	where $d_{\rm H}$ denotes the Hausdorff distance between the sets.
\end{remark}

\begin{remark}\label{uniqnEquiv-A-B}
By Theorem {\rm \ref{thm:con}},
conditions (a) and (b) in Theorem {\rm \ref{thm:main theorem2}}
for the solutions satisfying properties \eqref{RegReflSol-Prop0}--\eqref{RegReflSol-Prop1-1-1}
in Definition {\rm \ref{admisSolnDef}} are equivalent.
\end{remark}

\begin{remark}\label{uniqnSolAreAdmissible}
	We note that, under either one of
	conditions (a) and (b) in Theorem {\rm \ref{thm:main theorem2}}, the solution is an
	admissible solution. Indeed, in both cases, the solution satisfies  properties
	\eqref{RegReflSol-Prop0}--\eqref{RegReflSol-Prop1-1-1} in Definition \ref{admisSolnDef}.
	If, in addition, condition (b)
	holds, then the solution is admissible, directly from Definition {\rm \ref{admisSolnDef}}.
	Remark {\rm \ref{uniqnEquiv-A-B}} shows the same for the case when condition $(a)$ holds.
\end{remark}

The proof of Theorem \ref{thm:main theorem2} is obtained by showing the following
proposition on the existence and uniqueness of a family of
admissible solutions that are continuous with respect to $\theta_{\rm w}$,
containing a given admissible solution.

\begin{proposition}\label{thm:main theorem}
	Fix $(\rho_0, \rho_1, \gamma)$. Define interval $I:=(\theta_{\rm w}^{\rm d},\frac{\pi}{2}]$ when $u_1\leq c_1$ and
	$I:=(\theta_{\rm w}^{\rm c},\frac{\pi}{2}]$ when $u_1>c_1$.
	For every admissible solution $\varphi^*$ with a wedge angle $\theta_{\rm w}^*\in I$,
	there exists a family
	$$
	\mathfrak{S}=\{(\varphi, \theta_{\rm w})\;:\;\; \theta_{\rm w}\in I, \;
	\varphi\in C^{0,1}(\Lambda(\theta_{\rm w}))\}
	$$
	such that
	\begin{equation}\label{startFromSoluiton}
		(\varphi^*, \theta_{\rm w}^*)\in \mathfrak{S},
	\end{equation}
	and $\mathfrak{S}$ satisfies the following properties{\rm :}
	
	\smallskip	
	\begin{enumerate}[\rm (a)]
		\item \label{thm:main theorem-i1}
		For each $\theta_{\rm w}\in I$, there exists one and only one
		pair $(\varphi, \theta_{\rm w})\in\mathfrak{S}$.
		Then we can define $\varphi^{(\theta_{\rm w})}:=\varphi$
		if $(\varphi, \theta_{\rm w})\in\mathfrak{S}$.
		
		\smallskip
		\item Each $\varphi^{(\theta_{\rm w})}$ is an admissible
		solution corresponding to the wedge angle $\theta_{\rm w}$.
		
		\smallskip
		\item \label{thm:main theorem-i3}
		$\varphi^{(\frac{\pi}{2})}$ is the normal shock reflection solution $($\emph{see \S 3.1 in \cite{cf-annofmath20101067regularreflection} for the definition}$)$.
		
		\item \label{thm:main theorem-i4}
		$\varphi^{(\theta_{\rm w})}$ is continuous with respect to the wedge angle $\theta_{\rm w}\in I$
		in the $C^1$--norm as in Remark {\rm \ref{distBetweenRegReflSol}}.
	\end{enumerate}
	Moreover, a family $\mathfrak{S}$
	satisfying properties \eqref{thm:main theorem-i1}--\eqref{thm:main theorem-i4} listed above
	$($but without requiring \eqref{startFromSoluiton}$)$
	is unique. That is, if there are two families $\mathfrak{S}_1$ and $\mathfrak{S}_2$ satisfying
	properties \eqref{thm:main theorem-i1}--\eqref{thm:main theorem-i4},
	then $\mathfrak{S}_1=\mathfrak{S}_2$.
	Thus, the family $\mathfrak{S}$ contains all the admissible solutions for
	all $\theta_{\rm w}\in I$.
\end{proposition}

Proposition \ref{thm:main theorem} directly implies Theorem \ref{thm:main theorem2}.

\smallskip
Proposition \ref{thm:main theorem} is proved by showing the local uniqueness and existence
of admissible solutions.

As we have discussed in the introduction, the outline of the uniqueness proof ({\it i.e.} Proposition  \ref{thm:main theorem})
is the following:
If there are two different admissible solutions, defined by the potential functions
$\varphi$ and $\hat\varphi$,
for some wedge angle
$\theta_{\rm w}^*\in I\setminus\{ \frac\pi 2\}$, it suffices to:

\smallskip
\begin{itemize}
	\item[\rm (i)] construct  continuous  families
	of solutions parametrized by the wedge angle $\theta_{\rm w} \in [\theta_{\rm w}^*, \frac\pi 2]$,
	starting from $\varphi$ and $\hat\varphi$, respectively, in the norm discussed
	in Remark \ref{distBetweenRegReflSol};
	
	\smallskip
	\item[\rm (ii)] prove {\it local uniqueness}: If two admissible
	solutions with the same wedge angle
	are close in the norm given in the second line of (\ref{distBetweenRegReflSol-Eqn}), then they are equal.
\end{itemize}
Combining this with the fact that, by Remark \ref{admisEquiv1-2}, both families converge to the normal
reflection as $\theta_{\rm w}\to \frac\pi 2-$ , we obtain a contradiction; see more
details in \S \ref{UniquenessProofSubseq} below.
Furthermore, the continuous family defined above can be extended to all
$\theta_{\rm w}\in I$, hence
determining the family $\mathfrak{S}$ in Proposition \ref{thm:main theorem}.

In order to construct the continuous family of solutions
$\mathfrak{S}$
described in
Proposition \ref{thm:main theorem}, starting from the given solution
$\varphi=\varphi^{(\theta_{\rm w}^*)}$,
it suffices to show that any given admissible
solution can be perturbed, that is, an admissible solution can be constructed to be
close to  $\varphi$ for all wedge angles
that are sufficiently close to $\theta_{\rm w}^*$.
More precisely, using the mapping of admissible
solutions to the functions on the unit square discussed in Remark \ref{distBetweenRegReflSol},
we work in an appropriately weighted and scaled $C^{2,\alpha}$ space on $Q^{\rm iter}$.
We choose this function space according to the norms and the other quantities
in the {\it a priori} estimates for admissible solutions in \cite{cf-book2014shockreflection},
mapped to $Q^{\rm iter}$.
Denote the norm in this space by $\|\cdot \|^*$. Thus, we consider space $C_*(Q^{\rm iter})$, which
is completion of $C^\infty(\overline{Q^{\rm iter}})$ with respect to norm $\|\cdot\|^*$.
This space satisfies
$$
C_*(Q^{\rm iter})\subset C^{1,\alpha}(\overline{Q^{\rm iter}})\cap C^{2,\alpha}(Q^{\rm iter}).
$$
For any admissible solution $\varphi$, the corresponding function $u$ on $Q^{\rm iter}$
satisfies $u\in C_*(Q^{\rm iter})$.
Now we state the local existence assertion.

\begin{proposition}[Local existence]\label{localPerturbExist}
	Fix any admissible
	solution $(\hat\varphi, \hat\theta_{\rm w})$ with $\hat\theta_{\rm w}\in I$.
	Then, for every sufficiently small $\varepsilon>0$, there is
	$\delta>0$ with the following property{\rm :} For
	each $\theta_{\rm w}\in [\hat\theta_{\rm w}-\delta, \hat\theta_{\rm w}+\delta)]\cap I$,
	there exists an admissible solution $\varphi$ such that $u$ and $\hat u$ on
	$Q^{\rm iter}$ corresponding to $\varphi$ and $\hat\varphi$, respectively, satisfy
	$$
	\|u-\hat u\|^*<\varepsilon.
	$$
\end{proposition}
Note that, if $\varepsilon$ is sufficiently small, the solutions obtained
in Proposition \ref{localPerturbExist}
are unique for each wedge angle, by the local uniqueness.

\smallskip
Thus, to prove Proposition  \ref{thm:main theorem}, it suffices to prove the local uniqueness,
as well as the local existence in the sense of Proposition \ref{localPerturbExist}. See
\S \ref{UniquenessProofSubseq}
below for more details in the proof of Proposition \ref{localPerturbExist} from
these properties. In fact, from
Remark \ref{eqnInOmega-rmk}, we study these questions for the free boundary problem
\eqref{equ:potential flow equation} and \eqref{equ:boundary-RH-RegRefl}--\eqref{equ:boundary-ofOmega-RegRefl},
where the unknowns are $\varphi$ in $\Omega$ and $\Gsh$. Moreover,
the admissible solutions satisfy property
{\rm (\ref{RegReflSol-Prop1})} in Definition \ref{admisSolnDef}, from which equation
\eqref{equ:potential flow equation} is strictly elliptic
in $\overline\Omega\setminus\overline{\Gso}$ in the supersonic and
sonic cases $|D\varphi_2(P_0)|\ge c_2$ and uniformly elliptic in $\overline\Omega$
in the subsonic case $|D\varphi_2(P_0)|< c_2$.

\smallskip
The proofs of the local existence and uniqueness are different for the following two cases:

\smallskip
\begin{enumerate}[(a)]
	\item \label{supesonocCase}
	Supersonic and subsonic-near-sonic case{\rm :} $|D\varphi_2(P_0)|>(1-\sigma)c_2$,
	
	\smallskip
	\item \label{subsonocCase}
	Subsonic-away-from-sonic case{\rm :}  $|D\varphi_2(P_0)|\le (1-\sigma)c_2$,
\end{enumerate}
where $\sigma>0$ depends on $(\rho_0,\rho_1,\gamma)$ and is such that,
for the wedge angles satisfying $(1-\sigma)c_2\le |D\varphi_2(P_0)|\le 1$ (which are
the subsonic-near-sonic and sonic cases), the admissible solutions are $C^{2,\alpha}$ up to $P_0$
according to \cite{cf-book2014shockreflection}.

The reason for the different proofs for cases \eqref{supesonocCase} and \eqref{subsonocCase}
is that, in the supersonic and sonic case, the degenerate ellipticity
of equation \eqref{equ:potential flow equation} near $\Gso$ or $P_0$
makes it difficult to use the linearization of  problem
\eqref{equ:potential flow equation} and \eqref{equ:boundary-RH-RegRefl}--\eqref{equ:boundary-ofOmega-RegRefl}
for the application of the implicit function theorem
which would imply both the local existence and uniqueness.
On the other hand, under the conditions stated in case \eqref{supesonocCase},
$\varphi$ is $C^{1,1}$ up to $\Gso$ in the supersonic case (by
Theorem \ref{ExistRegRefl-Th1}(\ref{RegReflSol-Prop0})), and $C^2$ up to $P_0$
in the subsonic-near-sonic and sonic cases; this higher regularity allows us
to use the different methods described below.
In the subsonic-away-from-sonic case (\ref{subsonocCase}),
the known regularity up to $P_0$ is $C^{1,\alpha}$, {\it i.e.}
lower than that in case \eqref{supesonocCase},
but equation \eqref{equ:potential flow equation} is uniformly elliptic
in $\overline\Omega$; this allows to analyze the linearization of problem
\eqref{equ:potential flow equation} and \eqref{equ:boundary-RH-RegRefl}--\eqref{equ:boundary-ofOmega-RegRefl}
at $\varphi$, and thus obtain the local uniqueness
and existence by the implicit function theorem.

It remains to discuss the proof of the local uniqueness and existence
in the supersonic and subsonic-near-sonic case \eqref{supesonocCase}.
The outline of this proof is in \S 5.1--\S 5.2 below.

\subsection{Local uniqueness in the supersonic and subsonic-near-sonic case \eqref{supesonocCase}}
Assume that $\varphi$ and $\varphi^*$ are regular shock
reflection solutions for the same wedge angle $\theta_{\rm w}$,
which are $C^{1,1}$ up to $\overline\Gso$ (where we denote $\overline\Gso=\{P_0\}$
in the subsonic and sonic cases) and satisfy the properties listed in
Theorem \ref{thm:main theorem2}. Let
${\Omega}$ and ${\Omega^*}$ be respectively their elliptic regions, and let
$\Gsh$ and $\Gsh^*$ be respectively their reflected shocks.
We recall that $\varphi$ and $\varphi^*$ satisfy \eqref{equ:potential flow equation} and
\eqref{equ:boundary-RH-RegRefl}--\eqref{equ:boundary-ofOmega-RegRefl}
in $\Omega$ and $\Omega^*$, respectively.

Let $\hat{\Omega}:=\Omega\cap\Omega^*$, and let $\hGsh:=\partial\hat{\Omega}\cap(\Gsh^*\cup\Gsh)$.
We now show that, under the following assumption:
\begin{equation}\label{assunique}
	\|\varphi-\varphi^*\|_{C^1(\hat{\Omega})}+\|\varphi-\varphi_1\|_{C^0((\Omega\cup\Omega^*)\backslash\hat{\Omega})}
	+\|\varphi^*-\varphi_1\|_{C^0((\Omega\cup\Omega^*)\backslash\hat{\Omega})}\leq\delta_2,
\end{equation}
the function, $\delta\varphi:=\varphi-\varphi^*$, satisfies
the boundary condition:
\begin{align}\label{Bou:Lin Shock}
	\mathcal{M}(\delta\varphi)=\beta_{\bn}(\delta\varphi)_{\bn} +
	\beta_{\bt}(\delta\varphi)_{\bt}+\vartheta\delta\varphi=0 \qquad\mbox{on the inner shock $\hGsh$},
\end{align}
with
\begin{align} \label{Bou:Lin Shock-prop}
	\beta_{\bn}>0, \qquad \vartheta<0,
\end{align}
where $\bn$ is the unit inner normal and $\bt$ is the unit tangent on $\hGsh$.
We note that the property that $\vartheta<0$ in (\ref{Bou:Lin Shock-prop})
is obtained by using the convexity of $\Gsh^*$ and $\Gsh$.

Also, it follows from \cite{cf-comparison principles} that $\delta\varphi$
satisfies a homogeneous linear
elliptic equation in $\hat\Omega$
for which the comparison principles hold.
Properties (\ref{Bou:Lin Shock-prop}), combined with methods of
\cite{cf-comparison principles}, show that Hopf's lemma holds for
$\delta\varphi$ on $\hGsh$.
Finally, $\delta\varphi$ satisfies the
homogeneous Neumann condition on $(\partial\hat\Omega\cap \partial\Lambda)\setminus \{P_3\}$,
and $\delta\varphi=0$ on $\Gso$.

These facts ensure that $\delta\varphi\equiv 0$ in $\hat \Omega$.
From this, we can show
\begin{equation}\label{uniqueness}
	\Omega^*=\Omega,\qquad\quad\;  \varphi=\varphi^*\quad\mbox{in }\Omega.
\end{equation}
This completes the proof of the local uniqueness.

\begin{remark}
	We remark that, due to the issue that the regularity of
	$\varphi$ at the reflection point $P_0$ is only $C^{1,\alpha}$ for the
	subsonic-away-sonic reflection case $|D\varphi_2(P_0)|\le(1-\sigma)c_2$,
	we cannot apply
	this argument.
	However, as we discussed earlier, the implicit function theorem can be applied in that case.
\end{remark}

\subsection{Local existence in the supersonic and subsonic-near-sonic case \eqref{supesonocCase}}

Now we discuss the proof of the local existence, Proposition \ref{localPerturbExist}.
The existence
of a solution is obtained by the application of the Leray-Schauder degree
theory \cite[\S 13.6($A4^*$)]{zeidler-book1986nonlinearanalysisFixedpoint}; see also
\cite[\S 3.4]{cf-book2014shockreflection}.

In order to apply the degree theory, the iteration set should be bounded
and open in an appropriate function space (in fact, in its
product with the parameter space, {\it i.e.}
interval $[\hat\theta_{\rm w}-\delta, \hat\theta_{\rm w}+\delta]\cap I$ of the wedge angles),
the iteration map should be defined and continuous on the closure of the iteration
set, and any fixed point of the iteration map should not occur on the boundary
of the iteration set.
We choose this function space according to the norms and the other quantities
in the {\it a priori} estimates.
Moreover, since we have to use the same function space
for all values of the parameters,
and the functions require to have the same domain, we define the
iteration set in terms of the functions on the unit square $Q^{\rm iter}$,
which are related to the admissible solutions
by the mapping described in Remark \ref{distBetweenRegReflSol}.
The function space
is $C_*(Q^{\rm iter})$, introduced above.
Let $\hat u$ be the function on $Q^{\rm iter}$ corresponding
to the admissible solution $\hat\varphi$ for the wedge angle $\hat\theta_{\rm w}$ in
Proposition \ref{localPerturbExist}.
In order to prove the existence result in
Proposition \ref{localPerturbExist} for given $\varepsilon$ and $\delta$,
we define the iteration set by
\begin{equation}\label{iterSetDef}
	\mathcal{K}_{\varepsilon,\delta}^{(\hat u,\hat\theta_{\rm w})}:=
	\{(u, \theta_{\rm w})\in C_*(Q^{\rm iter})\times ([\hat\theta_{\rm w}-\delta,\; \hat\theta_{\rm w}+\delta]\cap I) \;\;:\;\;
	\|u-\hat u\|^*<\varepsilon\}.
\end{equation}
From its definition, the iteration set is non-empty, open (in the subspace topology) and bounded in
$C_*(Q^{\rm iter})\times ([\hat\theta_w-\delta,\; \hat\theta_w+\delta]\cap I)$.

We also define the iteration set for each wedge angle
$\theta_{\rm w}\in [\hat\theta_{\rm w}-\delta,\; \hat\theta_{\rm w}+\delta]\cap I$ by
\begin{equation}\label{iterSetDef-Thet}
	\mathcal{K}_{\varepsilon}^{(\hat u,\hat\theta_{\rm w})}(\theta_{\rm w}):=
	\{u\in C_*(Q^{\rm iter}) \;\;:\;\;(u, \theta_{\rm w})\in \mathcal{K}_{\varepsilon,\delta}^{(\hat u,\hat\theta_{\rm w})}\}.
\end{equation}
To prove Proposition \ref{localPerturbExist}, we need to show the existence of an admissible solution in
$\mathcal{K}_{\varepsilon,\theta_{\rm w}}^{(\hat u,\hat\theta_{\rm w})}$ for each $\theta_{\rm w}\in
[\hat\theta_{\rm w}-\delta,\; \hat\theta_{\rm w}+\delta]\cap I$ if $\varepsilon$ is small,
depending on $(\rho_0, \rho_1, \gamma, \hat\theta_{\rm w})$, and  $\delta$ is small, depending on
$\varepsilon$ and $(\rho_0, \rho_1, \gamma, \hat\theta_{\rm w})$.

The iteration map $\IterMap$ is defined as follows:

Given  $(u, \theta_{\rm w})\in \overline\IterSet$,
define the corresponding {\it elliptic domain}
$\Omega=\Omega(u, \theta_{\rm w})$ by mapping
from the unit square $\ItRg$ to
the {\it physical plane}, as discussed in Remark \ref{distBetweenRegReflSol}.
This determines iteration $\Gsh$ and function $\varphi$ in $\Omega$,
depending on $(u, \theta_{\rm w})$.
We set up a boundary
value problem in $\Omega$ for a {\it new iteration potential} $\tilde\varphi$
by modifying problem \eqref{equ:potential flow equation} and
\eqref{equ:boundary-RH-RegRefl}--\eqref{equ:boundary-ofOmega-RegRefl},
by partially substituting $\varphi$ into the coefficients
of \eqref{equ:potential flow equation}, and making other modifications
including the ellipticity cutoff in the equation.

In the supersonic and sonic cases,
the modified equation is elliptic in $\overline{\Omega}\setminus\Gso$,
degenerate near $\Gso$ (or $P_0$ in the sonic case), and nonlinear near  $\Gso$.
In the subsonic case, the modified equation is linear
and uniformly elliptic in $\overline\Omega$.

In all the supersonic, sonic, and subsonic cases, we prescribe one condition on $\Gsh$,
which is an oblique derivative condition, by combining the two conditions in
\eqref{equ:boundary-RH-RegRefl} and partially substituting $\varphi$ into the coefficients of
the main terms.

Let $\tilde\varphi$  be the solution of the boundary
value problem in $\Omega$. We show that $\tilde\varphi$ gains
the regularity in comparison with $\varphi$.
Then we define $\tilde u$ on $\ItRg$ by mapping $\tilde \varphi$ back
in such a way that the gain-in-regularity of the solution is preserved,
which is needed in order to have the compactness of the iteration map.
This requires some care, since the original mapping between $\ItRg$
and the {\it physical domain} is defined by $u$ and hence
has a lower regularity.
Then the iteration map is defined by
$$
\IterMap(u, \theta_{\rm w})=\tilde u.
$$
The boundary value problem in the definition of $\IterMap$
is defined so that, at the fixed point $u=\tilde u$, its solution satisfies
the potential flow equation \eqref{equ:potential flow equation}
with the ellipticity cutoff in a small neighborhood
of $\Gso$ in the supersonic case, both the Rankine-Hugoniot conditions
\eqref{equ:boundary-RH-RegRefl} on $\Gsh$,
and the boundary condition \eqref{equ:boundary-ofOmega-RegRefl} on $\Gw\cup\Symm$.
On the sonic arc $\Gso$ in the supersonic case and at $P_0$ in the subsonic
and sonic cases, we need two conditions: $\tilde\varphi=\varphi_2$ and $D\tilde\varphi=D\varphi_2$.
However, we can prescribe only one condition on the fixed boundary. We choose
the Dirichlet condition $\tilde\varphi=\varphi_2$ on  $\Gso$ in the supersonic case and at $P_0$ in the subsonic
and sonic cases,
and prove that $D\tilde{\varphi}=D\varphi_2$ on $\Gso$ or at $P_0$ holds for the solution
of the iteration problem  for the fixed point.

Then we prove the following facts:

\smallskip
(i) Any fixed point $u = \IterMap(u, \theta_{\rm w})$, mapped to the {\it physical plane},
is an admissible solution $\varphi$.
For that, we remove the ellipticity cutoff
and prove the inequalities and monotonicity properties in the definition of the admissible solutions
for the regions
and the wedge angles where they are not readily known from the definition of the iteration set.

(ii) The iteration map is continuous on $\overline\IterSet$
and compact. We prove this by using the gain-in-regularity of the solution of
the iteration boundary value problem.

\smallskip
(iii) Any fixed point of the iteration map cannot occur on the boundary of
the iteration set if $\delta$ is small depending on $\varepsilon$ and $(\rho_0, \rho_1, \gamma)$.
Now we discuss this step in more details:

\smallskip
The {\it small} iteration set (\ref{iterSetDef-Thet})
is the first key difference between this proof of the local existence and the proof of the existence
of admissible solutions in \cite{cf-book2014shockreflection}, which is also obtained by
the Leray-Schauder degree argument. In \cite{cf-book2014shockreflection}, the continuity of
admissible solutions with respect to $\theta_{\rm w}$ was not studied; for this reason, the iteration
set is chosen to be {\it large} for the wedge angles away from $\frac\pi 2$.
That is,
the iteration set for such a wedge angle is defined by the bounds in
the appropriate norms related to the {\it a priori} estimates
and by the lower bounds of certain directional derivatives,
corresponding to the strict monotonicity properties
so that the actual solution cannot be on the boundary of the iteration set
according to the {\it a priori} estimates.  In the present case of
{\it small} iteration set (\ref{iterSetDef-Thet}), a different approach is developed,
based on the local uniqueness and compactness of admissible solutions
shown in \cite{cf-book2014shockreflection}.
That is, fixing small $\varepsilon>0$, and assuming that, for any $\delta>0$,
there exists an admissible
solution $\tilde\varphi$ for the wedge angle $\tilde\theta_{\rm w}$ such that
$|\tilde\theta_{\rm w}-\hat\theta_{\rm w}|\le\delta$ and
$\|\tilde u-\hat u\|^*=\varepsilon$,
we obtain a sequence of admissible solutions and their
wedge angles
$(\varphi^{(i)}, \theta_{\rm w}^{(i)})$ with
$\theta_{\rm w}^{(i)}\to \hat\theta_{\rm w}$ and $\|u^{(i)}-\hat u\|^*=\varepsilon$.
Then, using the compactness of admissible solutions, we can send to a limit for a subsequence
so that an admissible solution $\bar\varphi$ is obtained for the wedge angle $\hat\theta_{\rm w}$ such that
$\|\bar u-\hat u\|^*=\varepsilon$. This contradicts the local uniqueness if $\varepsilon$ is small.

Now the Leray-Schauder
degree theory  guarantees that
the fixed point index:
\begin{equation}\label{fixedPtIndex}
	{\rm Ind}(\IterMap(\cdot, \theta_{\rm w}), \;
	\overline{\mathcal{K}_{\varepsilon}^{(\hat u,\hat\theta_{\rm w})}(\theta_{\rm w})})
\end{equation}
of the iteration map on
the iteration set (for given $\theta_{\rm w}$)
is independent of
the wedge angle $\theta_{\rm w}\in [\hat\theta_{\rm w}-\delta,\; \hat\theta_{\rm w}+\delta]\cap I$.

It remains to show that, at some wedge angle, index (\ref{fixedPtIndex}) is non-zero. We
show that, for the wedge angle $\hat\theta_{\rm w}$,
\begin{equation*}
	{\rm Ind}(\IterMap(\cdot, \hat \theta_{\rm w}), \;
	\overline{\mathcal{K}_{\varepsilon}^{(\hat u,\hat{\theta}_{\rm w})}(\hat\theta_{\rm w})})=1.
\end{equation*}
We prove this by showing that
\begin{equation}\label{itermap-is-const}
	\IterMap(v, \hat \theta_{\rm w})=\hat u
	\qquad\mbox{for each $v\in \mathcal{K}_{\varepsilon}^{(\hat{u},\hat{\theta}_{\rm w})}(\hat\theta_{\rm w})$}.
\end{equation}
This means that the iteration boundary value problem in domain $\Omega(v, \hat\theta_{\rm w})$
defined by every $v\in \mathcal{K}_{\varepsilon}^{(\hat{u},\hat{\theta}_{\rm w})}(\hat\theta_{\rm w})$
has the unique solution $\hat\varphi$ (in fact, its carefully defined
extension from $\Omega(\hat u, \hat\theta_{\rm w})$).
This step is another key difference from the existence proof of
admissible solutions in \cite{cf-book2014shockreflection}.
In \cite{cf-book2014shockreflection},
the iteration set includes the normal reflection $\varphi^{\rm normal}$ for $\theta_{\rm w}=\frac\pi 2$,
and property (\ref{itermap-is-const}) is shown for $\theta_{\rm w}=\frac\pi 2$ and $u^{\rm normal}$ on the
right-hand side. Since $\varphi^{\rm normal}$ is an explicitly known uniform state, globally defined,
showing (\ref{itermap-is-const}) is straightforward for the normal reflection,
and does not require defining its extension,
or any special properties of the coefficients of the iteration problem.
In the present case,
when $\hat\varphi$ is an arbitrary admissible solution,
this step is much more involved, and requires
an extension of $\hat\varphi$ from $\Omega$ to a larger region (so that the extension satisfies certain properties)
and some careful definition of the coefficients of the iteration equation and the boundary condition on $\Gsh$,
for which we need at least the $C^{1,1}$--regularity of $\varphi$ near $\Gso$.
Thus, our method works for the supersonic and subsonic-near-sonic case;
however, it does not readily work for the subsonic-away-from-sonic case
(for this reason, in this case, we use a different approach
as we discussed above).

This completes the proof of the local existence of supersonic
and subsonic-near-sonic reflection solutions.

\subsection{Proof of Proposition \ref{thm:main theorem}}
\label{UniquenessProofSubseq}

Based on the local uniqueness and existence, we employ the compactness of admissible solutions
proved in \cite{cf-book2014shockreflection} to
conclude that, for every
admissible solution $\varphi^*$ with the wedge angle $\theta_{\rm w}^*\in I$,
a family $\mathfrak{S}$
with the properties listed in Proposition  \ref{thm:main theorem} exists.

It remains to prove the uniqueness of admissible solutions for each wedge angle.

For a given wedge angle $\theta_{\rm w}$ as in Theorem \ref{thm:main theorem2},
assume that there
are two admissible solutions $\varphi$ and $\tilde\varphi$ corresponding to the wedge angle
$\theta_{\rm w}^*$.
Let $\mathfrak{S}$ and $\tilde{\mathfrak{S}}$ be the continuous
families with $(\varphi,\theta_w^*)\in\mathfrak{S}$
and $(\tilde{\varphi},\theta_w^*)\in \tilde{\mathfrak{S}}$
in Proposition  \ref{thm:main theorem}.
Let $\mathfrak{A}$ be the set of
all $\theta_{\rm w}\in [\theta_{\rm w}^*, \frac\pi 2]$ such that $\varphi^{\theta_{\rm w}}=\tilde\varphi^{\theta_{\rm w}}$.
Since $\frac\pi 2\in \mathfrak{A}$
by \eqref{thm:main theorem-i3} of Proposition \ref{thm:main theorem}, it follows that $ \mathfrak{A}\ne\emptyset$.
The continuity of both families $\mathfrak{S}$ and $\tilde{\mathfrak{S}}$ with respect to $\theta_w$
implies that
$ \mathfrak{A}$ is closed. Also, by the assumption above, $\theta_{\rm w}^*\notin \mathfrak{A}$.
Denote $\theta_{\rm w}^{\rm inf}:=\inf \mathfrak{A}$, then $\theta_{\rm w}^{\rm inf}\in (\theta_{\rm w}^*, \frac\pi 2]$.
Now, using the continuity of families $\mathfrak{S}$ and $\tilde{\mathfrak{S}}$, we can show
that, choosing $\theta_{\rm w}\in (\theta_{\rm w}^*, \theta_{\rm w}^{\rm inf})$ to be sufficiently close
to $\theta_{\rm w}^{\rm inf}$, we obtain that $\varphi^{(\theta_{\rm w})}=\tilde\varphi^{(\theta_{\rm w})}$ by the
local uniqueness property.
This contradicts the definition of $\theta_{\rm w}^{\rm inf}$.


\end{document}